\documentclass[11pt]{article}
\usepackage{amsfonts}
\usepackage{graphics}
\usepackage{eucal}
\usepackage{graphics}

\newtheorem{Proposition}{Proposition }
\newtheorem{Lemma}{Lemma}

\newcommand{\newsection}[1]{\section{#1}\setcounter{equation}{0}}

\newcommand{\for}[1]{(\ref{#1})}

\newcommand{\bib}[1]{\textup{\cite{#1}}}

\begin{document}

\newcommand\ugarr{\!\!\!&=&\!\!\!}   

\renewcommand{\a}{\alfa}
\renewcommand{\b}{\beta}
\newcommand{\e}{\epsilon}
\newcommand{\g}{\gamma}
\newcommand{\p}{\varphi}
\renewcommand{\P}{\Phi}
\newcommand{\s}{\sigma}
\renewcommand{\a}{\alpha}
\renewcommand{\d}{\delta}
\renewcommand{\o}{\omega}
\renewcommand{\r}{\rho}
\renewcommand{\l}{\lambda}
\renewcommand{\L}{\Lambda}
\renewcommand{\O}{\Omega}


\newcommand{\beq}[1]{\begin{equation}\label{#1}}
\newcommand{\eeq}{\end{equation}}
\newcommand{\Prop}[1]{\acapon\begin{Proposition}\label{#1}}
\newcommand{\eProp}{\end{Proposition}}
\newcommand{\Lem}[1]{\acapon\begin{Lemma}\label{#1}}
\newcommand{\eLem}{\end{Lemma}}

\newcommand{\cA}{\mathcal{A}}                 
\newcommand{\cB}{\mathcal{B}}
\newcommand{\cC}{\mathcal{C}}                 
\newcommand{\cF}{\mathcal{F}}                 
\newcommand{\cG}{\mathcal{G}}                 
\newcommand{\cH}{\mathcal{H}}
\newcommand{\cI}{\mathcal{I}}                 
\newcommand{\cL}{\mathcal{L}}                 
\newcommand{\cM}{\mathcal{M}}                 
\newcommand{\cO}{\mathcal{O}}                 
\newcommand{\cP}{\mathcal{P}}                 
\newcommand{\cR}{\mathcal{R}}                 
\newcommand{\cS}{\mathcal{S}}                 
\newcommand{\cT}{\mathcal{T}}                 
\newcommand{\cV}{\mathcal{V}}                 
\newcommand{\cW}{\mathcal{W}}                 
\newcommand{\cX}{\mathcal{X}}

\newcommand{\deq}{$\interi_4\per\interi_2$--equivalent\ }
\newcommand{\deqp}{$\interi_4\per\interi_2$--equivalent}
\newcommand{\deqence}{$\interi_4\per\interi_2$--equivalence\ }
\newcommand{\deqencep}{$\interi_4\per\interi_2$--equivalence}
\newcommand{\naturali}{\mathbb{N}}
\newcommand{\Sp}[1]{ \mathrm{Sp}(#1)}
\newcommand{\fou}[2]{\langle #1\rangle_{#2}}

\marginparsep=2truemm
\marginparwidth=1.5truecm
\newcommand{\fcomment}[1]{\marginpar{\scriptsize\tt #1 (F)}}
\newcommand{\dcomment}[1]{\marginpar{\scriptsize\sl #1 (D)}}


\newcommand{\acapo}{\vspace{2ex}}
\newcommand{\acapon}{\vspace{2ex}\noindent}
\newcommand{\qqquad}{\quad\quad\quad}

\newcommand{\ins}[3]{\vbox to0pt{\kern-#2 \hbox{\kern#1
#3}\vss}\nointerlineskip}    

\newcommand{\toro}{\mathbb{T}}
\newcommand{\complessi}{\mathbb{C}}
\newcommand{\reali}{\mathbb{R}}
\newcommand{\interi}{\mathbb{Z}}

\newcommand{\rdue}{\reali^2}
\newcommand{\rtre}{\reali^3}
\newcommand{\rsei}{(\reali^3)^2}

\newcommand{\unity}{\mathbf{1}}


\newcommand{\ug}{\;=\;}
\newcommand{\piu}{\;+\;}		
\newcommand\meno{\;-\;}	
\newcommand\TO{\;\to\;}	
\newcommand\GE{\;\ge\;}	
\newcommand\LE{\;\le\;}	

\newcommand{\per}{\times}		
\renewcommand{\bar}{\overline}		

\newcommand{\step}[1]{^{(#1)}}
\newcommand{\der}[2]{\frac{\partial#1}{\partial#2}}
\newcommand{\dder}[3]{\frac{\partial^2#1}{\partial#2\partial#3}}

\newcommand{\grad}{\mathrm{grad}} 
\newcommand{\const}{\mathrm{const}} 
\newcommand\rank{\mathrm{rank\,}}
\newcommand\diag{\mathrm{diag\,}}
\newcommand\tr{\mathrm{tr\,}}
\renewcommand\skew{\mathrm{skew\,}}
\newcommand\sign{\mathrm{sign\,}}

\newcommand\fineprova{\ \vrule height6pt width6pt depth0pt}

\newcommand{\SKEWTRE}{\mathrm{so(3)}}
\newcommand{\OTRE}{\mathrm{O(3)}}
\newcommand{\LTRE}{\mathrm{L(3)}}
\newcommand{\GLTRE}{\mathrm{GL(3)}}
\newcommand{\SLTRE}{\mathrm{SL(3)}}
\newcommand{\sltre}{\mathrm{sl(3)}}
\newcommand{\SLDTRE}{\mathrm{SLD(3)}}
\newcommand{\sldtre}{\mathrm{sld(3)}}
\newcommand{\SOTRE}{\mathrm{SO(3)}}
\newcommand{\sotre}{\mathrm{so(3)}}
\newcommand{\PS}[2]{\langle#1,#2\rangle}
\newcommand{\bigPS}[2]{\big\langle#1,#2\big\rangle}
\newcommand{\BigPS}[2]{\Big\langle#1,#2\Big\rangle}
\newcommand{\colvect}[2]{\left(\matrix{#1\cr#2\cr}\right)}

\newcommand{\sub}[1]{_{(#1)}}


\newcommand{\D}{D}       
\newcommand{\MR}{G^R}      
\newcommand{\MS}{G^S}      
\newcommand{\wMS}{G}       
\newcommand{\wMR}{G}       
\newcommand{\GG}{\Psi}     

\newcommand{\Sdue}{$\mathrm{S}_2$ }
\newcommand{\Stre}{$\mathrm{S}_3$ }
\newcommand{\Sduep}{$\mathrm{S}_2$}
\newcommand{\Strep}{$\mathrm{S}_3$}

\newcommand{\cBSdue}{\cB_{\mathrm{S}_2}}
\newcommand{\cBStre}{\cB_{\mathrm{S}_3}}
\newcommand{\cBI}{\cB_{\mathrm{I}}}
\newcommand{\cBII}{\cB_{\mathrm{II}}}
\newcommand{\cBIII}{\cB_{\mathrm{III}}}

\newcommand{\twopig}{}
\newcommand{\pig}{{1\over2}}


\newcommand{\mytable}[4]{
\begin{table}
\center{\small
\begin{tabular}{||l|l||l|l||l|l||l|l||l|l||} 
\hline\hline 
\multicolumn{10}{||c||}{}\\
\multicolumn{10}{||c||}{{\em #2 }} \\
\multicolumn{10}{||c||}{}\\
\multicolumn{1}{||c|}{#3} & \multicolumn{1}{|c||}{#4} & 
\multicolumn{1}{c|}{#3} & \multicolumn{1}{|c||}{#4} &
\multicolumn{1}{c|}{#3} & \multicolumn{1}{|c||}{#4} & 
\multicolumn{1}{c|}{#3} & \multicolumn{1}{|c||}{#4} & 
\multicolumn{1}{c|}{#3} & \multicolumn{1}{|c||}{#4} \\ 
\hline
\input{#1}
\hline\hline 
\end{tabular}
}
\end{table}
}

\newcommand{\I}{{\mathbb I}}
\newcommand{\E}{{\mathbb E}}
\newcommand{\F}{{\mathbb F}}
\newcommand{\K}{{\mathbb K}}
\newcommand{\J}{{\mathbb J}}

\def\mommap{\mathrm{\bf J}}


\newcommand{\norm}[1]{\left \vert #1 \right \vert}
\newcommand{\proj}[2]{{\mathbb P}_{#1} #2}
\newcommand{\lp}{\left (}
\newcommand{\rp}{\right )}
\newcommand{\la}{\left \langle}
\newcommand{\ra}{\right \rangle}
\newcommand{\lsb}{\left [}
\newcommand{\rsb}{\right ]}
\newcommand{\lcb}{\left \{}
\newcommand{\rcb}{\right \}}
\newcommand{\R}{{\mathbb R}}

\newcommand{\sands}{\qquad \mbox{and} \qquad}
\newcommand{\half}{{\textstyle {1 \over 2}}}
\newcommand{\smallfrac}[2]{{\textstyle {#1 \over #2}}}

\newcommand{\ddt}{\smallfrac{d\ }{dt}}
\newcommand{\dds}{\smallfrac{d\ }{ds}}
\newcommand{\ddeps}{\smallfrac{d\ }{d\e}}
\newcommand{\ddlambda}{\smallfrac{d\ }{d\l}}


\newcommand\bList{
\begin{list}{}{\leftmargin2em\labelwidth1em\labelsep.5em\itemindent0em
\topsep0ex\itemsep-.8ex} }
\newcommand\eList{\end{list}}

\newcommand\bListList{
\begin{list}{}{\leftmargin2.5em\labelwidth1.7em\labelsep.5em\itemindent0em
\topsep0ex\itemsep-.8ex} }

\newcommand\finelemma{\raisebox{4pt}
           {\framebox[6pt][5pt]{\hbox to 6pt{}}
}}

\def\uncatcodespecials{\def\do##1{\catcode`##1=12 }\dospecials}
\def\sic{\begingroup\tt\uncatcodespecials\obeylines
  \obeyspaces\doverbatim}
\newcount\balance
{\catcode`<=1 \catcode`>=2 \catcode`\{=12 \catcode`\}=12
  \long\gdef\doverbatim{<\balance=1\verbatimloop>      
  \long\gdef\verbatimloop#1<\def\next<#1\verbatimloop>%
    \if#1{\advance\balance by1
    \else\if#1}\advance\balance by-1
     \ifnum\balance=0\let\next=\endgroup\fi\fi\fi\next>>
%

\newcommand{\ml}{m_l}
\newcommand{\mr}{m_r}


\title{\bf Stability Properties of the Riemann Ellipsoids }

\author{
Francesco Fass\`o\thanks{Universit\`a di Padova, Dipartimento di
Matematica Pura e Applicata, Via G. Belzoni 7, 35131 Padova,
Italy. (E--mail: fasso@math.unipd.it). Partially supported by
a CNR-NATO Fellowship, by DOE contract DEFG03-95ER25245-A000
while visiting the Mathematics Department, University of
California at Santa Cruz, and by MURST protocol 9802261238\_008.}
\and
Debra Lewis\thanks{Mathematics Department, University of
California, Santa Cruz, CA 95064, USA. (E--mail:
lewis@math.ucsc.edu). Partially supported by NSF grant 
DMS-9802378 and the Santa Fe Institute. }
}

\vskip 1truecm
\date{\small(January 31, 2000)}

\maketitle

\vskip 1truecm
\begin{abstract} \noindent\small
We study the ellipticity and the ``Nekhoroshev stability"
(stability properties for finite, but very long, time scales) of
the Riemann ellipsoids. We provide numerical evidence that the
regions of ellipticity of the ellipsoids of types II and III are
larger than those found by Chandrasekhar in the 60's and that
all Riemann ellipsoids, except a finite number of codimension
one subfamilies, are Nekhoroshev--stable. We base our analysis
on a Hamiltonian formulation of the problem on a covering space,
using recent results from Hamiltonian perturbation theory.

\vskip2truecm
\end{abstract}

\vskip2truecm

\newsection{Introduction} 

{\bf A. } The Dirichlet problem consists of the ``linear''
motions of an ideal, incompressible, homogeneous,
self--gravitating fluid mass having an ellipsoidal shape at any
instant. This problem originated in attempts to determine the
shapes of the planets, and has a long and distinguished history
(see e.g. chapter~1 of \bib{chandrasekhar69} and references
therein). Let us only remark here that Newton and Maclaurin
established the existence of rigidly rotating oblate spheroids
and that Jacobi described rigidly rotating asymmetric
ellipsoids. Dirichlet~\bib{dirichlet860} proved the existence of
``linear'' fluid motions, such that the position $x(t,y)$ at
time $t$ of any material particle $y$ is given by
\beq{motion}
    x(t,y) = F(t) y \,,\qquad F(t)\in\SLTRE \,.
\eeq
If, as we assume, the reference configuration of the fluid mass
is a ball of radius $\rho$, then the free surface of the fluid
determined by the motion \for{motion} is at every instant an
ellipsoid with semiaxes $\rho a_1$, $\rho a_2$, $\rho a_3$,
where $a_1$, $a_2$ and $a_3$ are the singular values of the
matrix $F$ (namely, the eigenvalues of $\sqrt{FF^T}$).

The Maclaurin spheroids and Jacobi ellipsoids are steady
rigid motions of the Dirichlet problem. Dedekind identified a new
class of steady, but non--rigid, motions of the problem; they
were constructed by transposition of the Jacobi ellipsoids, thus
bringing the particle relabelling symmetry of the system into
play. Riemann reformulated the equations of motion of the
Dirichlet problem in a particularly convenient form for studying
asymmetric configurations (as we shall show, this
reformulation corresponds to the passage to a four--to--one
covering space). Riemann also classified all steady motions of
the system, showing that they are generated by combinations of
spatial and body symmetries, and studied their Lyapunov (or
nonlinear) orbital stability.

The steady motions found by Riemann, which are called {\em
Riemann ellipsoids,} are motions \for{motion} such that
$$
    F(t) \ug \exp(t\O_l)\, A \,\exp(-t\O_r) 
$$
for some constant diagonal matrix $A=\diag(a_1,a_2,a_3)$ and
constant antisymmetric matrices $\O_l$ and~$\O_r$. Except for
the special case of the Jacobi ellipsoids, for which $\O_r=0$,
Riemann ellipsoids do not perform rigid motions: the motion is
the composition of an internal rotation, a stretch along the
principal axes, and a spatial rotation; the free surface retains
a rotating ellipsoidal shape, but the fluid particles describe
rosette--shaped motions which are either closed (periodic) or
open (quasi--periodic) depending on the two angular frequencies
$\o_l$ and $\o_r$ (that is, the vectors corresponding to the
antisymmetric matrices $\O_l$ and~$\O_r$). Riemann proved that
there are five types of Riemann ellipsoids with distinct
semiaxes, characterized by different relations between the angular
velocities and the semiaxes. For two types (hereafter referred to
as \Sdue and \Strep), both $\o_l$ and $\o_r$ are parallel to the
same principal axis of the ellipsoid; this axis is either the
shortest (\Strep) or the middle (\Sduep) one. For the other
three types, called I, II, and III, both $\o_l$ and $\o_r$
belong to one of the two principal planes containing the longest
axis of the ellipsoid.\footnote{The distinction between \Sdue
and \Stre ellipsoids that we make here is not standard;
Chandrasekhar calls them collectively ellipsoids of type S.}

Since the Riemann ellipsoids are reduced equilibria, namely equilibria 
modulo symmetries, of the Dirichlet problem, it is meaningful to 
study their stability modulo symmetries, which we will
frequently refer to simply as stability. 
Riemann found that all ellipsoids of type \Stre and a
subfamily of the \Sdue ellipsoids are stable. Riemann's analysis
does not lead to any definite conclusion about the (in)stability of
the remaining ellipsoids, because they are saddle points of his
Lyapunov function.

More recently, Chandrasekhar studied the `ellipticity', or
`spectral stability', of the Riemann ellipsoids. 
We say here that a (relative) equilibrium is elliptic if the 
eigenvalues of the linearization at the equilibrium are purely
imaginary; note that, accordingly, we treat the case of zero 
eigenvalues as elliptic, even though this is not customary. 
Chandrasekhar
\bib{chandrasekhar65}\bib{chandrasekhar66}\bib{chandrasekhar69}
numerically found a region of ellipticity of the ellipsoids of
type \Sdue larger than the region of Lyapunov stability found by
Riemann. He also found nonempty regions of elliptic ellipsoids
of type I and III, but none of type II. As we shall see, some of
these conclusions are not entirely correct, and in fact require
substantial modification.

\acapo

\acapon
{\bf B. } For a Hamiltonian system such as the Dirichlet problem,
ellipticity is a necessary, but not sufficient, condition for
Lyapunov stability. Therefore, the ellipticity analysis is of 
twofold interest: on one hand, the absence of ellipticity
implies Lyapunov instability; on the other hand, it opens up the
problem of the stability properties of those Riemann ellipsoids
that are elliptic, but of undetermined Lyapunov stability.
The present article is mainly devoted to a study of the latter 
question. 

First of all, as a prerequisite for this analysis, we repeat the
Lyapunov stability and the ellipticity analysis of the Riemann
ellipsoids. While our analysis confirms Riemann's Lyapunov
stability results, we find significant discrepancies with
Chandrasekhar's results: there is a nonempty region of elliptic
ellipsoids of type II and a region of elliptic ellipsoids of
type III that is substantially larger than the one reported by
Chandrasekhar. Therefore, {\em the regions of known (Lyapunov)
instability of the ellipsoids of types II and III are
substantially smaller than those found by Chandrasekhar.}
Moreover, the ellipticity regions of the ellipsoids of type I,
II and III have a finer and more complex structure than was
previously realized.

The main goal of the article is to use the methods of
Hamiltonian perturbation theory to investigate the long--time
stability properties of the elliptic Riemann ellipsoids of
undetermined Lyapunov stability. In fact, the Dirichlet problem
is a Hamiltonian system with symmetry and the Riemann ellipsoids
are the relative equilibria of this system. They correspond to
the equilibria of the reduced system, which has four degrees of
freedom for generic values of the momenta and three degrees of
freedom in special cases. (The reduced phase space is in fact an
orbifold, with a singular set containing the Riemann ellipsoids
of type S, but this difficulty is overcome by passing to Riemann's
covering space.) Since the reduced system has
more than two degrees of freedom, KAM theory cannot be used to
show the stability of its equilibria. Therefore, we shall base
our investigation on Nekhoroshev theory, which allows one to
obtain ``practical'' stability results, namely, to prove that
motions starting sufficiently near the equilibrium remain near
it for times growing exponentially with (some power) of the
initial distance; specifically, we shall use some recent results
obtained in \bib{benfasguz98} (see also \bib{fasguzben98},
\bib{niederman98}, \bib{guzfasben98}). The results of our
investigation indicate that {\em all elliptic Riemann
ellipsoids are Nekhoroshev--stable, with the possible exception
of a finite number of codimension one families of ellipsoids 
satisfying certain low--order resonance conditions.} 

In view of the complexity of the system, we resort to numerical
methods at certain steps in the ellipticity analysis and in the
construction of the normal forms for the reduced Hamiltonian. 
We numerically compute the eigenvalues of the
linearization and numerically construct the normal forms for a
finite (though rather large and, we think, significant) number
of values of the momentum. Therefore, the forementioned results
are not rigorous, but we think that our analysis provides strong
evidence for our conclusions.

\acapo

\acapon
{\bf C. } Our analysis is based on a Hamiltonian formulation (on
a covering manifold) of the Dirichlet problem, which is briefly
developed in Section~2. For the sake of brevity, we skip most of
the computational details and focus on the conceptual features.
The Riemann ellipsoids are described in Section 3. Overall,
Sections 2 and 3 contain a treatment of the Dirichlet problem 
which, although concise, is complete and detailed; this seems
necessary, because the existing literature does not seem to
include a comparably comprehensive treatment.

Sections 4--6 contain the ellipticity and the
Nekhoroshev--stability analyses. In Section~4 we review the
results on Lyapunov stability and instability and describe the
results of our ellipticity analysis; a detailed comparison with
Chandrasekhar's results is included in Appendix~B. We briefly
review the necessary notions for the Nekhoroshev--stability
analysis in Section~5 and describe the results for the Riemann
ellipsoids in Section~6. A Conclusion follows.

All of the numerical computations have been performed using the
software package {\sl Mathematica}.

\acapo

\newsection{The Dirichlet Problem}

{\bf A. Dirichlet's equations. } Dirichlet \bib{dirichlet860}
showed that the linear motions \for{motion} form an invariant
submanifold of the phase space of an ideal, incompressible,
homogeneous, self--gravitating fluid, with the boundary
condition that the pressure is constant at the free surface. He
expressed the restriction of the hydrodynamic equations of motion to
this submanifold as equations for the curve $F(t)\in\SLTRE$,
obtaining a second order system on $\SLTRE$. 

To describe these results, we introduce the following notation,
which will be used throughout the article. We denote by $\cdot$
and $\|\ \|$ the Euclidean scalar product and norm on $\rtre$ and
by $D_\r$ the ball $\{y\in\reali^3:\|y\|\le\r\}$. We use the
standard inner product and norm on $L(3)$ given by
\beq{PS}
    \PS A B \ug \tr(A B^T) \quad \mathrm{and} \quad 
    \big|A\big|^2 \ug \PS{A}{A} \,.
\eeq
(Note that the subspaces of symmetric and antisymmetric matrices
are $\PS{\;}{\,}$--orthogonal complements of one another in $L(3)$.) 
We shall frequently identify $\rtre$ and $\SKEWTRE$ by means of the
isomorphism $\o \mapsto \hat \o$, where $\hat \o \, \xi= \o \times
\xi$ for all $\xi\in\rtre$. Note that 
$\PS {\hat\o}{\hat\xi} \ug 2\, \o\cdot\xi$
for all $\o,\,\xi\in\rtre$. 
Following Riemann, we utilize the ``singular value decomposition'' 
of matrices, nowadays a standard tool in numerical linear algebra: If
$F\in\SLTRE$ and $a_1,a_2,a_3$ are its singular values, then
there exist matrices $U_l, U_r\in\SOTRE$ such that
$$
   F = U_l A U_r^T  \,, \qquad A\ug \diag(a_1,a_2,a_3) \,.
$$
We shall refer to the triplet $(A,U_l,U_r)$ as a ``SVD'' of
the matrix $F$. The matrices $A$, $U_l$ and $U_r$ are not
unique, but one can clearly fix $A$ by ordering its entries,
e.g. $a_1\ge a_2\ge a_3$ (see also Proposition \ref{covering} below).
Finally, we set
$$
    \cV(a) \ug -2\pi g \int_0^\infty 
    \Big[(s+a^2_1)(s+a^2_2)(s+a^2_3)\Big]^{-1/2} ds  \,,
$$
where $g$ denotes the gravitational constant; by a slight abuse
of notation, we shall regard $\cV$ and its derivatives as
functions either of $a=(a_1,a_2,a_3)$ or of the diagonal matrix
$A=\diag(a_1,a_2,a_3)$.

Dirichlet's result can now be expressed as follows: {\em
Consider a curve $t\mapsto F(t)\subset \SLTRE$ defined on some
open interval $\cI\subset\reali$ and let $(A(t),U_l(t),U_r(t))$
be any singular value decomposition of $F(t)$. Then the map
$$
    x\,:\, \cI\per D_\r\;\to\; \rtre \,,\quad  
    (t,y)\;\mapsto\; F(t)y 
$$
is a solution of the hydrodynamical equation for an ideal,
incompressible, homogeneous, self--gravitating fluid with
constant pressure at the boundary iff
\beq{eq-motion-1} 
   \proj{F}\Big[
   \ddot F + \twopig U_l\cV'(A)U_r^T
   \Big]= 0 \,,
\eeq
where 
$\cV'=\diag\Big(\der{\cV}{a_1},\der{\cV}{a_2},\der{\cV}{a_3}\Big)$
and $\proj{F}{G} = G - {1 \over 3} \PS G F F^{-T}$ for any $G\in
L(3)$. } Note that the matrix $U_l\cV'(A)U_r^T$ is independent 
of the particular SVD of $F$ used; this is proven by observing
that, if $U,V\in\SOTRE$ and $A=UA'V$ with $A$ and $A'$ diagonal,
then $\cV'(A)=\cV'(UA'V)=U\cV'(A')V$.
In fact, Dirichlet used the polar decomposition, not the SVD,
when writing the force term in the equations of motion.

\acapo 
{\em Remark: } In his work on Dirichlet problem, Chandrasekhar
required that the relative pressure $p-p_o$, where $p_o$ is the
constant pressure on the surface of the fluid, be nonnegative at
any point of the fluid. Following Riemann, we do not make this
assumption. (There are many fluids that support even very high
tensile states.)

\acapo
Since the gravitational forces are conservative and invariant
under spatial rotations and particle relabelling, equations
\for{eq-motion-1} inherit the conservation of energy, angular
momentum and circulation from the hydrodynamic equations. For
future reference, we observe that if we set
$c={m\r^2\over5}$, where $m$ denotes the mass of the fluid, then
the energy of any motion
\for{motion} is $c\, \big[{1\over2}\big|\dot F\big|^2 + \twopig
\cV(a)\big]$ and the angular momentum is the vector $k$ satisfying
$\hat k \ug 2\,c\,\, \mbox{skew} \big(\dot FF^T \big)$; here
$\mbox{skew}(G)={1\over2}(G-G^T)$. The conservation of
circulation is equivalent to the constancy of $\mbox{skew}
\big(F^T \dot F \big)$.

\acapo

\acapon
{\bf B. Riemann's equation. } Shortly after Dirichlet's work
appeared, Riemann \bib{riemann860} rewrote equation
\for{eq-motion-1} using what is nowadays called the singular
value decomposition of matrices and left--trivializing the factors
$T\SOTRE$, namely expressing the tangent vectors $\dot U_l$ and 
$\dot U_r$ in terms of the antisymmetric matrices $\O_l=U_l^T
\dot U_l$ and $\O_r = U_r^T \dot U_r$. Once these substitutions 
have been made, equation \for{eq-motion-1} takes the form
\beq{eq-motion-2}
  {\mathbb P}_A \left[ \ddot A 
   \piu 2(\O_l\dot A -\dot A\O_r) + (\dot{\O}_l A - A\dot{\O}_r)
   \piu (\O^2_l A - 2\O_l A \O_r + A \O^2_r)  
    + \twopig \cV'(A) \right] = 0 
\eeq
(which should be compared with equations ($\a$) in
\bib{riemann860}). 
This equation determines a
second order system on $\cA\per\SOTRE\per\SOTRE$, where
$$
  \cA \ug
  \Big\{ \diag\Big(a_1,a_2,(a_1a_2)^{-1}\Big) \,:\;
         a_1>a_2>a_1^{-1/2} \Big\}
$$
(if $A$ has repeated eigenvalues, $\dot\O_l$ and $\dot\O_r$ are
not uniquely determined by \for{eq-motion-2}). This system is
equivalent to the restriction of Dirichlet's equation
\for{eq-motion-1} to the submanifold
$$
   Q = \{F\in \SLTRE: F\; \mathrm{has\ distinct\ singular\
        values}\} \,.
$$
More precisely, because of the nonuniqueness of the SVD, the two 
systems are related not by a diffeomorphism, but by a 
four--to--one covering map:

\Prop{covering} 
The map
$$
   \cC: \cA\times\SOTRE\times\SOTRE \to Q \,,\qquad
        (A,U_l,U_r)\mapsto U_l\, A \, U_r^T 
$$
is a 4--fold covering map which maps solutions of equation
\for{eq-motion-2} into solutions of equation \for{eq-motion-1}.
\eProp

\acapon
{\bf Proof. } If $(A,U_l,U_r)$
and $(A,V_l,V_r)$ are any two SVDs of $F\in Q$, with
$A=\diag(a_1,a_2,a_3)$, $a_1>a_2>a_3$, then
$V_l \ug R_i U_l$ and $V_r \ug R_i U_r$
for some $i=0,\ldots,3$, where 
\beq{cover_matrices}
R_0=\unity, \quad R_1=\diag(1,-1,-1), \quad
R_2=\diag(-1,1,-1), \quad \mbox{and} \quad R_3=\diag(-1,-1,1).
\eeq
Thus, the map $\cC$ is four--to--one. That it is a submersion is
proven, for instance, by showing that the determinant of its
Jacobian is nonzero. (Using $a_1,a_2$ as coordinates on $\cA$
and left--trivializing the factors $T\SOTRE$, one obtains the 
determinant $(a_1^2-a_2^2)(a_1^2-a_3^2)(a_2^2-a_3^2)$.) \fineprova

\acapo
In Riemann's formulation, the use of the SVD has the advantage
that the absence of the matrices $U_l$ and $U_r$ from the
equation \for{eq-motion-2} makes manifest the $\SOTRE\per\SOTRE$
invariance of the Dirichlet problem. From our perspective, the
formulation on the covering also has the advantage that, once we
pass to the Hamiltonian formulation, the symplectically reduced
phase space associated to the cotangent bundle of the covering
is a manifold, while reduction of the system on $T^*Q$ produces
singularities.

\acapo
{\em Remarks: } (i) The submanifold $Q$ is not invariant under
the dynamics of the Dirichlet problem (axisymmetric
configurations can evolve into asymmetric ones, and vice versa)
but this is not a problem for our analysis: the Riemann
ellipsoids with asymmetric configurations clearly belong to $Q$
for all time and, since $Q$ is open, it will be a consequence of
our stability estimates that nearby motions also remain in $Q$
for the time under consideration.

(ii) The restriction to the case of distinct singular
values excludes all axisymmetric Riemann ellipsoids. The stability 
properties of the Maclaurin spheroids, which rotate rigidly 
about the axis of symmetry, are completely known (see
\bib{riemann860}\bib{chandrasekhar69}\bib{lewis93}); the remaining
axisymmetric ellipsoids merit further investigation.

\acapon

\acapon 
{\bf C. Hamiltonian formulation. } We now describe the
Hamiltonian formulation of Riemann's equation \for{eq-motion-2}
on the cotangent bundle of the covering manifold
$\cA\per\SOTRE\per\SOTRE$. We shall use the first two diagonal
entries $a_1$ and $a_2$ as coordinates on the set $\cA$, by
means of the diffeomorphism $\a:\cB\to\cA$, where
$$
  \cB \ug 
  \Big\{(b_1,b_2)\in\rdue:\, b_1>b_2>b_1^{-1/2} \Big\} 
  \qquad \mbox{and} \qquad
   \a(b_1,b_2) := \diag\Big(b_1,b_2,{1\over b_1 b_2} \Big) \,.
$$
(We denote by $b=(b_1,b_2)$ the first two singular values, so as
to distinguish $b\in\rdue$ from $a\in\rtre$; wherever we write
$b_3$, it stands for $(b_1b_2)^{-1}$, not for an independent
variable.) Furthermore we identify the cotangent bundle of each
factor $\SOTRE$ with $\SOTRE\per\rtre$ as follows: we first
identify the cotangent space $T^*_U\SOTRE$ with $T_U\SOTRE$ by
means of the inner product ${1\over2}\PS{\;}{\;}$ (see \for{PS})
and then identify $T_U \SOTRE$ with $\rtre$ by means of the $\
\hat{}\ $ isomorphism, thus associating to each $\mu\in
T^*_U\SOTRE$ the unique $m\in\rtre$ satisfying
$$
   \mu \cdot U \hat \o \ug {1\over2}\,\PS{\hat m}{\hat\o} 
                \ug m\cdot\o \qquad \forall \ \o\in\rtre\,.
$$
After these identifications, we work on the 
sixteen--dimensional manifold
$$
   \cM \ug \cB \per\rdue \per(\SOTRE)^2\per(\rtre)^2 \,,
$$
which is diffeomorphic to $T^*\big(\cA \times \SOTRE \times
\SOTRE\big)$. We shall denote by $(b,c,U,m)$ the elements of
$\cM$, with $b\in\cB$, $c\in\rdue$, $U=(U_l,U_r)\in(\SOTRE^2)$, and
$m=(m_l,m_r)\in(\rtre)^2$ and use the induced scalar and vector
products on $(\rtre)^2$ given by
$$
   (\o_l,\o_r)\cdot(\xi_l,\xi_r) \ug 
   \o_l\cdot\xi_l + \o_r\cdot\xi_r 
   \,,\qquad
   (\o_l,\o_r)\times(\xi_l,\xi_r) \ug
   (\o_l\times\xi_l,\o_r\times\xi_r) \,.
$$
When writing tangent vectors to $\cM$ we shall again
left--trivialize and identify the tangent spaces to the two
factors $\SOTRE$ with $\rtre$; thus, $(v,w,\o,n)$ denotes the
vector $(v,w,(U_l\hat\o_l,U_r\hat\o_r),n) \in T_{(b,c,U,m)}\cM$.
Finally, in order to avoid the introduction of a new symbol, we
write $\cV(b)$ to denote $\cV\big(b_1,b_2,{1\over b_1
b_2}\big)$.

\Prop{Ham} 
Riemann's equation \for{eq-motion-2} on $\cA\per\SOTRE\per\SOTRE$
is equivalent to the Hamiltonian system on the manifold $\cM$ 
defined by the Hamiltonian
\beq{Hamiltonian}
   H(b,c,U,m) \ug  
   \half \, c\cdot\K(b)c \piu \half \, m\cdot\J(b)m 
   + \twopig\cV(b) \,,
\eeq
where 
\begin{eqnarray*}
   \K(b) \ugarr 
   {1\over b_2^2 b_3^2 +b_1^2 b_3^2 +b_1^2b_2^2}
   \left( \matrix{ b_1^2 \lp b_2^2 + b_3^2 \rp &  -b_3 \cr
                  - b_3 & b_2^2 \lp b_1^2 + b_3^2 \rp  \cr }\right) \\
   \J(b) \ugarr
   \left(\matrix{ J_1(b) & J_2(b) \cr
                  J_2(b) & J_1(b) \cr}\right) 
   \quad \mathrm{with\ }\ 
   \cases{
      \!\!\!\!\!\!\!&
        $J_1(b)= 
        \diag\Big({b_2^2+b_3^2\over (b_2^2-b_3^2)^2}, 
                  {b_1^2+b_3^2\over(b_1^2-b_3^2)^2}, 
                  {b_1^2+b_2^2\over (b_1^2-b_2^2)^2} \Big)$ \cr
      \!\!\!\!\!\!\!&
        $J_2(b) =
        \diag\Big({2 b_2b_3 \over (b_2^2-b_3^2)^2}, 
                  {2 b_1b_3 \over(b_1^2-b_3^2)^2}, 
                  {2b_1b_2\over (b_1^2-b_2^2)^2} \Big)$ \,, \cr}       
\end{eqnarray*}
and the ``left--trivialized'' (in the sense just specified)
symplectic structure $\s$ is defined by
\beq{sigma}
   \s\big(b,c,U,m\big) 
    \big((v,w,\o,n), (v',w',\o',n')\big)     
   \ug v\cdot w' \meno v'\cdot w \piu
   \o\cdot n'\meno n\cdot \o' \piu  m\cdot(\o\per\o') \,.      
\eeq
\eProp

\acapo
This proposition can be proven by a direct computation showing 
that Hamilton's equations for the Hamiltonian system defined by 
\for{Hamiltonian} and \for{sigma} are equivalent to Riemann's equations.
Alternatively, one can observe that Riemann's equations are
equivalent to the Lagrangian system defined on $T\cA\per T\SOTRE
\per T\SOTRE$ by the Lagrangian 
\beq{Lagr-1}
   \cL((A,\dot A),U_l\hat\o_l,U_r \hat\o_r)\ug 
   \half \big|\dot A +\hat\o_l A-A\hat\o_r\big|^2 
   \meno \twopig \cV(A) 
\eeq
and construct the Hamiltonian system using the Legendre
transform. For the sake of brevity, we skip all details,
remarking only that the relation between `angular momenta' and
`angular velocities' is $\o=\J(b)m$.

\acapo
{\em Remark: }  The Lagrangian nature of Dirichlet's problem,
with Lagrangian equal to the difference of the kinetic and the
potential energies, is implicit in Riemann's work and has been
repeatedly used since then. One can verify by a direct
computation that the Euler--Lagrange equations for the
Lagrangian \for{Lagr-1} are equivalent to Riemann's equations
\for{eq-motion-2}. 

\acapo
The Hamiltonian system $(\cM,H,\s)$ of Proposition \ref{Ham} is
invariant under the symplectic action $\P$ of $\SOTRE\per\SOTRE$
on $\cM$ given by
$$
   \P_{(R_l,R_r)}\big(b,c,(U_l,U_r),m\big) \ug 
   \big(b,c,(R_lU_l,R_rU_r),m\big) \,.
$$
Since the group does not act on the factor $\cB\per\rdue$, 
while it acts on each factor $\SOTRE\per\rtre$ as the
left--trivialized cotangent lift of the left action of $\SOTRE$
on itself, all results from the elementary case of the left
action of $\SOTRE$ on $T^*\SOTRE$ apply (see, e.g., \bib{libe87}
or \bib{cushman+bates97}). Thus, one concludes that the action $\P$
has the momentum map
$$
   \mommap:\cM\to\rsei \,,\qquad 
   (a,c,U,m) \,\mapsto\,
   \big(U_lm_l,U_rm_r\big) \,.
$$
For any $\eta\in\rsei$, one may identify $\mommap^{-1}(\eta)$ and
$\cB\per\rdue\per(\SOTRE)^2$, in which case the immersion
$i_\eta:\mommap^{-1}(\eta) \to \cM$ is given by
$i_\eta\big(a,c,U)=\big(a,c,U,(U_l^T\eta_l, U_r^T\eta_r)\big)$.
The two components $U_lm_l$ and $U_rm_r$ of the momentum map
coincide, up to constant rescalings, with the total angular
momentum and the circulation; see Section 2.A.

The next Proposition, which describes the reduced system
$(P_\eta,\s_\eta,H_\eta)$, also follows immediately from standard
results for the reduction of an $SO(3)$ invariant Hamiltonian
system on $T^*\SOTRE$. We denote by $S^2_\r$ the sphere of
radius $\r$. As usual, $G_\eta$ denotes the isotropy subgroup of
$\eta$ with respect to the coadjoint action; in the case at
hand, $G_\eta = \{\mbox{rotations about $\eta_l$}\} \times
\{\mbox{rotations about $\eta_r$}\}$.

\Prop{action2} (i) For any $\eta=(\eta_l,\eta_r)\in\rsei$ with
$\eta_l\not=0\not=\eta_r$, the reduced phase space
$\mommap^{-1}(\eta)/G_\eta$ can be identified with
$$
   P_\eta\ug \cB\per\rdue\per \Big(S^2_{\|\eta_l\|}
                         \per S^2_{\|\eta_r\|}\Big) \,.
$$
After this identification, the canonical projection
$\pi_\eta:\mommap^{-1}(\eta)\to P_\eta$ is
$\pi_\eta\big(b,c,U\big) =
\big(b,c,(U_l^T\eta_l,U_r^T\eta_r)\big)$. The symplectic
structure $\s_\eta$ on $P_\eta$ defined by the equality
$i_\eta^*\s=\pi_\eta^*\s_\eta$ is
$$
  \s_\eta(b,c,m)\big((v,w,m\per\o),(v',w',m\per\o')\big) 
  \ug 
  v\cdot w' - v'\cdot w - m\cdot\o\per\o'  \,.
$$
The reduced Hamiltonian $H_\eta:P_\eta\to\reali$, defined by
$H_\eta\circ\pi_\eta=H\circ i_\eta$, is
\beq{Hmu}
    H_\eta(b,c,m) \ug \half \,c\cdot\K(b)c \piu
    \half m\cdot\J(b)m \piu \cV(b) \,.
\eeq

(ii) If $\eta=(\eta_l,0)$ or $\eta=(0,\eta_l)$ with $\eta_l
\not=0$, then the reduced phase space can be identified with
$\cB\per\rdue\per S^2_{\|\eta_l\|}$ and one has
$\pi_\eta(b,c,U_l)=(b,c,U_l^T\eta_l)$ and 
\begin{eqnarray*}
  \s_\eta(b,c,m_l)\big((v,w,m_l\per\o_l),(v',w',m_l\per\o'_l)\big)   \ug 
             v\cdot w' - v'\cdot w - m_l\cdot\o_l\per\o'_l  \\
  H_\eta(b,c,m_l) \ug \half \,c\cdot\K(b)c \piu
    \half \, m_l\cdot J_1(b)m_l \piu \cV(b) \,.
\end{eqnarray*}
\eProp

\acapo
The reduced system generically has four degrees of freedom, but
only three in the special case (ii), which arises in the case of
the so--called ``irrotational'' Riemann ellipsoids. The case
$\eta_l=\eta_r=0$ yields a reduced system with two degrees of
freedom, the only equilibrium of which is the stationary sphere.

In view of the application of Nekhoroshev theory, we note that
the reduced Hamiltonian $H_\eta$ is an analytic function. This
is seen by observing that the self--gravitational potential can
be written in the form
\[ 
   \cV(b) \ug -\,{4\,\pi\,g \over\sqrt{b_1^2-b_3^2}} \,
              F\left(\arccos\Big({b_3\over b_1}\Big) \,\Big|\, 
               \sqrt{b_1^2-b_2^2 \over b_1^2-b_3^2 }\, \right)
\] 
where $F(\varphi \, | \, k)$ is the incomplete elliptic integral
of the first kind, which is analytic for $\varphi \in \lp 0,
\frac{\pi}2 \rp$ and $k \in (0, 1)$ (see \bib{byrd54}, pg. 4--5,
299).

\acapon

\acapon 
{\bf D. Equivalence and symmetry of reduced systems. } As we now
discuss, (i) certain classes of reduced systems are equivalent and
(ii) the reduced systems have certain discrete symmetries. For the
sake of conciseness, we consider here only the generic case, in
which both $\eta_l$ and $\eta_r$ are nonzero. (The case in which
one of the momenta vanishes is recovered with obvious modifications.)
There are two distinct reasons for equivalence of different
reduced systems:
\bList
\item[E1.] The reduced system $(P_\eta,H_\eta,\s_\eta)$ depends on
$\eta=(\eta_l,\eta_r)$ only through the norms of $\eta_l$ 
and~$\eta_r$. This equivalence reflects the invariance of the 
Dirichlet problem under spatial rotation. 
\item[E2.] Dedekind proved that {\em if $F(t)$ is a solution of
the Dirichlet problem, then $F(t)^T$ is also a solution}
\bib{dedekind860}. (The proof follows immediately from the
observations that $\proj {F^T}{G^T} = (\proj F G)^T$ for any
matrix $G$ and that $U_l A U_r^T$ is a singular value
decomposition of $F$ iff $U_r A U_l^T$ is a singular value
decomposition of $F^T$.) Chandrasekhar calls `adjoint' the two
motions determined by $F(t)$ and $F(t)^T$. This invariance of
the Dirichlet problem has a counterpart in the Hamiltonian
formulation on the covering space $\cM$, which is invariant
under the symplectomorphism exchanging $U_l$ with $U_r$ 
and $m_l$ with $m_r$. This map takes solutions to solutions, but
exchanges the level sets of~$\mommap$. Correspondingly, the
reduced systems for the values $(\eta_l,\eta_r)$ and
$(\eta_r,\eta_l)$ of the momentum map are conjugated to one
another by the diffeomorphism $(b,c,(m_l,m_r)) \mapsto
(b,c,(m_r,m_l))$.
\end{list}

\noindent
Furthermore, each reduced system is invariant under two discrete
actions:
\bList
\item[i.] The $\interi_2$--action $(b,c,m) \mapsto (b,c,-m)$,
which reflects the invariance of the original system with
respect to the combination of spatial reflection and time
reversal. 
\item[ii.] The $\interi_4$--action $(b,c,m) \mapsto
(b,c,R_im)$, $i=0,\ldots,3$, where $R_0,\ldots,R_3$ are given by
\for{cover_matrices} and $R_i(m_l,m_r)=(R_im_l,R_im_r)$. The
points $(b,c,R_i m)$ are obtained from one another by a
reflection about one and the same coordinate axis on the two
spheres; see figure~\ref{Z4-Orbit}. Hence the $\interi_4$--orbit
of the point $(b,c,m)$ consists of two distinct points if $m$ is
parallel to a basis vector $(e_j,e_j)$, $j=1,2,3$, and of four
distinct points otherwise. 
\end{list}

\noindent
In the stability analysis, we need not distinguish equivalent
reduced systems, nor distinct points within the $\interi_4
\times \interi_2$ orbits, in each reduced system. A
$\interi_4\per\interi_2$--orbit consists of eight, four, or two
points, depending on the number of nonzero components of $m_l$
and $m_r$. We will say that all points in such an orbit are
\deq and regard two points in the reduced space as adjoint if
any two points in their $\interi_4\per\interi_2$--orbits are
adjoint.

By a {\em Riemann ellipsoid} we will mean an equivalence class,
under spatial rotations (E1) and transposition (E2), of motions of
the Dirichlet problem corresponding to an equilibrium of the reduced
system. Thus each Riemann ellipsoid can be identified with the
$\interi_4\per\interi_2$--orbit of an equilibrium of some reduced
system.\footnote{Note that, as was previously stated, we exclude
from our consideration all axisymmetric ellipsoids.}

\acapo
{\em Remark: } The $\interi_4$--symmetry is not present in
Dirichlet's system \for{eq-motion-1}, but is clearly inherited
from the passage to the fourfold covering manifold $\cM$ of
$T^*Q$. (The existence of this symmetry of Riemann's equation is
mentioned by Chandrasekhar (\bib{chandrasekhar69}, pg. 72), but
its origin is not identified there.) Since all points in a
$\interi_4$--orbit correspond to the same state of the fluid,
the passage to the covering introduces some redundancy, but
allows us to avoid singular reduction. In fact, Dirichlet's
equations \for{eq-motion-1} determine a
$\SOTRE\per\SOTRE$--invariant Hamiltonian system on $T^*Q$;
however, the action is not free and the reduced phase space is
the orbifold $\cB\per\rdue\per(S^2\per S^2)/\interi^4$, rather
than a smooth manifold. The singular set consists of all states
in which $m_l$ and $m_r$ are both parallel to the same axis of
the ellipsoid, and thus contains the Riemann ellipsoids of
type~S.
 
\begin{figure}
\vskip6truecm
\includegraphics{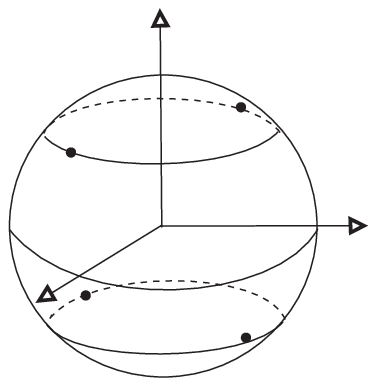}
\vskip-4truemm
\caption{\small Projection of a $\interi_4$--orbit on one of the spheres. 
\label{Z4-Orbit}
}
\vskip3truemm
\end{figure}

\acapon

\newsection{Riemann Ellipsoids}

We now study the equilibria of the reduced systems of
Proposition \ref{action2}, and thus the Riemann ellipsoids. The
results of this analysis are due to Riemann. As we mentioned in
the Introduction, there are five classes of Riemann ellipsoids
and the momenta $m_l$ and $m_r$ of the ellipsoids of each type
are completely determined by their semiaxes. Correspondingly, in
the set $\cB$ there are ``existence regions'' of the Riemann
ellipsoids of each type. In order to specify these regions, and
the momenta of the corresponding Riemann ellipsoids, we
introduce some notation. First we set
$$
   C_n(x,y,z) \ug 2\pi g \int_0^\infty 
    \Big[(s+x^2)(s+y^2)(s+z^2)\Big]^{-3/2} \, s^n ds
   \,,\qquad n=0,1,2 \,.
$$
The functions 
$$
    \MS_\pm(x,y,z) \ug {(x\mp z)^4\over xz} 
       \lsb \big[xy^2z \pm (x^2y^2-x^2z^2+y^2z^2)\big]C_1(x,y,z) \piu
        (xz\pm y^2)C_2(x,y,z) \rsb 
$$
determine the momenta of the ellipsoids of type S, while the functions
\begin{eqnarray*}
    \wMR(x,y,z) \ugarr 
        x^2\big(y^2-z^2\big)C_1(x,y,z) \piu 
        \big(y^2-4z^2\big)\big[z^2 C_1(x,y,z)+ C_2(x,y,z)\big]   \\
    \D(x,y,z) \ugarr x^2(y^2-z^2) + z^2(4z^2-y^2)                \\
    \MR_\pm(x,y,z) \ugarr 
        (y\mp z)^4 \lp x^2-(y\pm2z)^2 \rp \;
        {x^2-z^2 \over x^2-y^2} \;    
        {\wMR(x,y,z)\over \D(x,y,z)}
\end{eqnarray*}
determine the momenta of the ellipsoids of type I, II and III; we
refer to the latter three classes of ellipsoids collectively as type~R 
ellipsoids.\footnote{In Riemann's article, the study of the existence of 
these ellipsoids precedes that of the ellipsoids of type S.}
Using these functions we define the five sets
\begin{eqnarray}
   \cBSdue \ugarr \big\{b\in\cB\,:\, 
           G^S_-(b_1,b_2,b_3)\ge 0 \big\} \nonumber\\
   \cBStre \ugarr \big\{b\in\cB\,:\, 
           G^S_+(b_1,b_3,b_2)\ge 0 \big\} \nonumber\\
   \cBI \ugarr \big\{b\in\cB\,:\, 
           b_1\le 2b_2-b_3 \big\} \label{ER}\\
   \cBII \ugarr \big\{b\in\cB\,:\, 
           b_1\ge 2b_2+b_3 \,,\; 
           \D(b_1,b_3,b_2) <0   \big\}      \nonumber\\
   \cBIII \ugarr \big\{b\in\cB\,:\, 
           b_1\ge b_2+2b_3 \,,\; 
           \wMR(b_1,b_2,b_3)>0 \big\}       \nonumber
\end{eqnarray}
and five pairs of real--valued maps $\mu_\a^\pm$, one for each
type $\a$ of ellipsoids; these maps and their domains are
described in Table~1, where $\{e_1,e_2,e_3\}$ denotes the
standard basis vector of $\rtre$ and
\[
    N^\tau_\pm(x,y,z) \ug
    {\textstyle {1\over2}} \,\Big[\sqrt{G^\tau_+(x,y,z)} \,\pm\, 
                     \sqrt{G^\tau_-(x,y,z)} \, \Big]           
\qquad \qquad \tau = R, S \,.
\]

\vbox{
\center{
\begin{tabular}{||c|c|c||}
\hline\hline 
Type ($\a$) & Domain &$\mu_\a^\pm(b)$ \\ 
\hline\hline 
\Sdue
&$\cBSdue$ 
&$N^S_\pm(b_1,b_2,b_3)e_2$ 
\\ 
\hline 
\Stre
&$\cBStre$ 
&$N^S_\pm(b_1,b_3,b_2)e_3$ 
\\ 
\hline 
I
&$\cBI$ 
&$N^R_\pm(b_1,b_3,b_2)e_1 + N^R_\pm(b_3,b_1,b_2)e_3$ 
\\ 
\hline 
II
&$\cBII$
&$N^R_\pm(b_1,b_3,b_2)e_1 + N^R_\mp(b_3,b_1,b_2)e_3$
\\ 
\hline 
III
&$\cBIII$ 
&$N^R_\pm(b_1,b_2,b_3)e_1 + N^R_\mp(b_2,b_1,b_3)e_2$
\\ 
\hline\hline
\end{tabular}
}
\vskip3truemm
\centerline{\small\em Table 1: The momenta of the Riemann 
Ellipsoids}
}

\Prop{EqRel}

(i) A point $\big(b,c,m\big)\in P_\eta$ is an equilibrium of the
reduced system $(P_\eta,H_\eta,\s_\eta)$, $\eta_l \neq 0 \neq
\eta_r$, iff it is \deq to either the point
$\big(b,0,(\mu_\a^+(b),\mu_\a^-(b))\big)$ or the point
$\big(b,0,(\mu_\a^-(b),\mu_\a^+(b))\big)$ for some $b\in\cB_\a$
and some
$\a=\mathrm{S}_2,\mathrm{S}_3,\mathrm{I},\mathrm{II},\mathrm{III}$.

(ii) If one component of $\eta$ equals zero, then the equilibria
of the six--dimensional reduced system are the points \deq
to $(b,0,\mu_\a^+(b))$. This situation is met only
at the equilibria of type \Sdue satisfying $\MS_+(b_1,b_2,b_3)
= \MS_-(b_1,b_2,b_3)$ and those of type I satisfying
$b_1^2b_2^2+b_1^2b_3^2+b_2^2b_3^2-3b_2^4 =0$; in both cases,
$\mu_\a^+(b)\not=0=\mu_\a^-(b)$.
\eProp

\begin{figure}
\vskip17truecm
\includegraphics{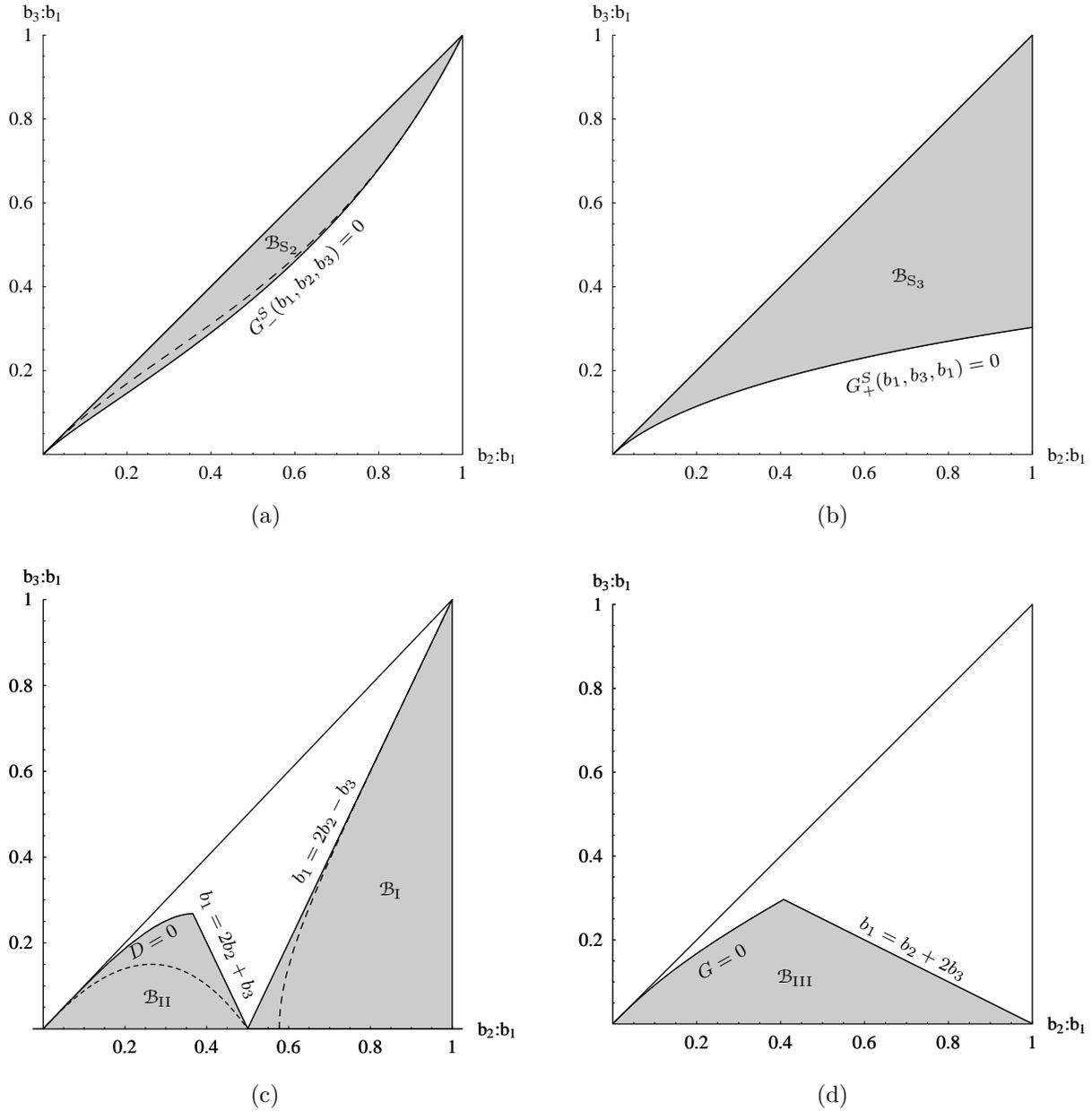}
{\scriptsize
\ins{37truemm}{88truemm}{\small (a)}
\ins{121truemm}{88truemm}{\small (b)}
\ins{39truemm}{128truemm}{$\cBSdue$}
\ins{36truemm}{131truemm}{\rotatebox{45}{$\MS_-(b_1,{b_2},{b_3})=0$}}
\ins{132truemm}{123truemm}{$\cBStre$}
\ins{125truemm}{111truemm}{\rotatebox{10}{$\MS_+(b_1,{b_3},{b_1})=0$}}
\ins{37truemm}{2truemm}{\small(c)}
\ins{21mm}{16mm}{$\cBII$}
\ins{29truemm}{31truemm}{\rotatebox{-65}{\scriptsize $b_1=2b_2+b_3$}}
\ins{18.5truemm}{26truemm}{\rotatebox{30}{\scriptsize $\D=0$}}
\ins{56mm}{32mm}{$\cBI$}
\ins{43mm}{48mm}{\rotatebox{63.435}{\scriptsize $b_1=2b_2-b_3$}}
\ins{121truemm}{2truemm}{\small(d)}
\ins{115mm}{19mm}{$\cBIII$}
\ins{127mm}{27mm}{\rotatebox{-26.6}{\scriptsize $b_1=b_2+2b_3$}}
\ins{103truemm}{23truemm}{\rotatebox{29.5}{\scriptsize $\wMR=0$}}
}

\vskip3truemm
\caption{\small The existence regions of the five types of
Riemann ellipsoids. The dashed curves in regions $\cBSdue$ and
$\cBI$ identify the irrotational ellipsoids and their adjoints
(for which $\nu_\a=0$); the former separates the coparallel
(lower) and the counterparallel (upper) \Sduep--ellipsoids. The
dashed curve in region $\cBII$ identifies the zero--pressure
ellipsoids.
\label{Exist}
}
\end{figure}

The proof of this proposition can be reconstructed from
Riemann's and Chandrasekhar's work. A sketch of the argument is
given in Appendix~A.

The existence regions \for{ER} are drawn in figures~\ref{Exist}
using the coordinates $(x,y)=({b_2/b_1},{b_3/b_1})$, which have
the advantage that the unbounded region $\cB$ is mapped onto the
bounded triangle $\{(x,y):\, 0<y<x<1\}$. The borders of the
various regions are determined by the definitions~\for{ER}. 
Note that, since the existence regions of ellipsoids of
different types have nonempty intersections, there exist Riemann
ellipsoids of up to four different types having equal semiaxes.

The Riemann ellipsoids belonging to a six--dimensional reduced
space have either zero angular momentum or zero circulation and
are called ``irrotational'' by Chandrasekhar. The existence of
irrotational ellipsoids of type \Sdue was known to Riemann, while
the existence of those of type~I seems not to have been noticed
before. Note that the curve of irrotational ellipsoids of type
\Sdue divides the region $\cBSdue$ into two parts, in one of which
the two momenta are coparallel, while in the other they are
counterparallel.\footnote{We say that two vectors $u$ and $v$
are coparallel (counterparallel) if $u\cdot v=\|u\|\,\|v\|$
($=-\|u\|\,\|v\|$).}

\acapo
{\em Remark: } Since Chandrasekhar requires that the relative
pressure be positive inside the ellipsoids, his region of
existence of the ellipsoids of type II is smaller than ours,
lying below the ``zero--pressure curve'' $\D(b_1,b_3,b_2)C_0 +
6b_2^2C_1 + 3C_2=0$ shown in figure~\ref{Exist}.c (see
\bib{chandrasekhar69}, chapter 7, formula (195)).

\acapo

\newsection{Stability and Ellipticity of the Riemann Ellipsoids}

{\bf A. Results. } We now review the known results about the
Lyapunov stability of the Riemann ellipsoids, regarded here as
equilibria in the reduced phase space, and we describe the 
ellipticity analysis. From now on, we restrict ourselves to the
generic case in which the reduced system has four degrees of 
freedom; the irrotational Riemann ellipsoids will be studied
elsewhere.

\Prop{LyapStab} 
(i) All \Strep--ellipsoids and all counterparallel
\Sduep--ellipsoids are Lyapunov stable. 

(ii) The coparallel \Sduep--ellipsoids for
$b\in\cBSdue\cap\mathrm{Int}(\cBIII)$ (i.e.,
$\wMR(b_1,b_2,b_3)>0$) are not elliptic and hence are Lyapunov
unstable.
\eProp

\acapo
The first statement is due to Riemann \bib{riemann860}, who
proved it using arguments essentially equivalent to using the
reduced Hamiltonian as a Lyapunov function. Riemann also showed
that the reduced Hamiltonian has a saddle point, and thus cannot
be used as a Lyapunov function, at any other ellipsoid. Riemann
interpreted this fact as an indication that these ellipsoids are
unstable (``labil"). However, as best we know this fact is
unproven, except in the case of the \Sdue ellipsoids considered
in statement (ii), the instability of which follows from
Chandrasekhar's work, since he showed that the linearization of
the equations of motion has an eigenvalue with nonzero real part
at these points; see section 4.B.

\acapo
{\em Remark: } Given the structure of the problem, the
Lyapunov stability of the equilibria of the reduced system
corresponds to orbital stability of the solutions of the
Dirichlet problem (see e.g. \bib{Patrick92}, \bib{Lewis92},
\bib{OrtegaRatiu98}).

\begin{figure}
\vskip7.2truecm
\includegraphics{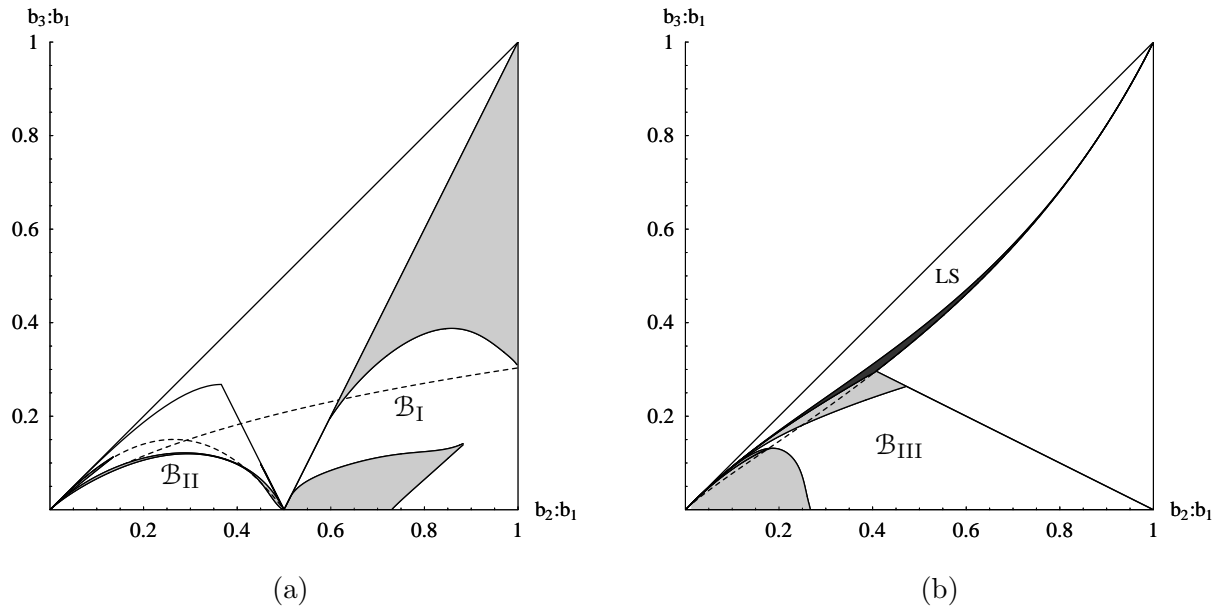}
\ins{22truemm}{17truemm}{$\mathrm{\cBII}$}
\ins{53truemm}{26truemm}{$\mathrm{\cBI}$}
\ins{117truemm}{21truemm}{$\mathrm{\cBIII}$}
\vskip-2truemm\centerline{(a) \hskip8cm (b)}
\vskip2truemm
\caption{\small Regions of ellipticity of the Riemann ellipsoids
of types I, II and III (light gray) and of the ellipsoids of
type \Sdue (dark gray). The dashed curve contained in region
$\cBII$ is the zero--pressure curve, while the one crossing
regions $\cBII$ and $\cBI$ is the border of the existence region
of the ellipsoids \Strep; the dashed curve in (b) is the border
of the existence region of the \Sduep.
The Riemann ellipsoids of types I, II and III in the unshaded
subregions of their existence regions are not elliptic and hence
are Lyapunov unstable. 
The counterparallel \Sdue ellipsoids, which lie in the unshaded
subregion of $\cBSdue$ marked LS, are Lyapunov stable; the
coparallel \Sdue ellipsoids in $\cBSdue\cap\mathrm{Int}(\cBIII)$ 
are Lyapunov unstable (in figure (b), this part of region 
$\cBSdue$ is `hidden' below $\cBIII$).
\label{Ell-1}
}
\end{figure}

\acapo
The ellipticity of the Riemann ellipsoids has been studied by
Chandrasekhar; most of his results are based on numeric 
calculations. We have
repeated the numerical analysis of the ellipticity, finding the
same results as Chandrasekhar for the \Sdue ellipsoids and
essentially the same results for the ellipsoids of type~I, but
significantly different results for the ellipsoids of types II
and III. The conclusions of our investigation are the following:

\acapon
{\bf Numerical Conclusions 1: } {\em (i) The coparallel
\Sduep--ellipsoids for $b$ in the complement of
$\cBSdue\cap\mathrm{Int}(\cBIII)$, i.e. $b$ satisfying
$\wMR(b_1,b_2,b_3)\le0$, are elliptic.

(ii) There are nonempty subregions of $\cBI$, $\cBII$ and
$\cBIII$ consisting of elliptic ellipsoids; these regions are
the shaded regions in figures~\ref{Ell-1} and \ref{Ell-2}. The
ellipsoids of types I, II and III in the complement of these
subregions are not elliptic and hence are Lyapunov unstable.
}

\acapo
This analysis has been carried out by analytically constructing
the linearization of the reduced Hamiltonian vector field at the
reduced equilibria and by numerically evaluating this matrix at
some tens of thousands of ellipsoids of each type. Details about
this analysis are given in the next subsection. (As explained
there, Riemann sketches an argument that would imply statement
(i).) Here we comment on these results.

\begin{figure}
\vskip18truecm
\includegraphics{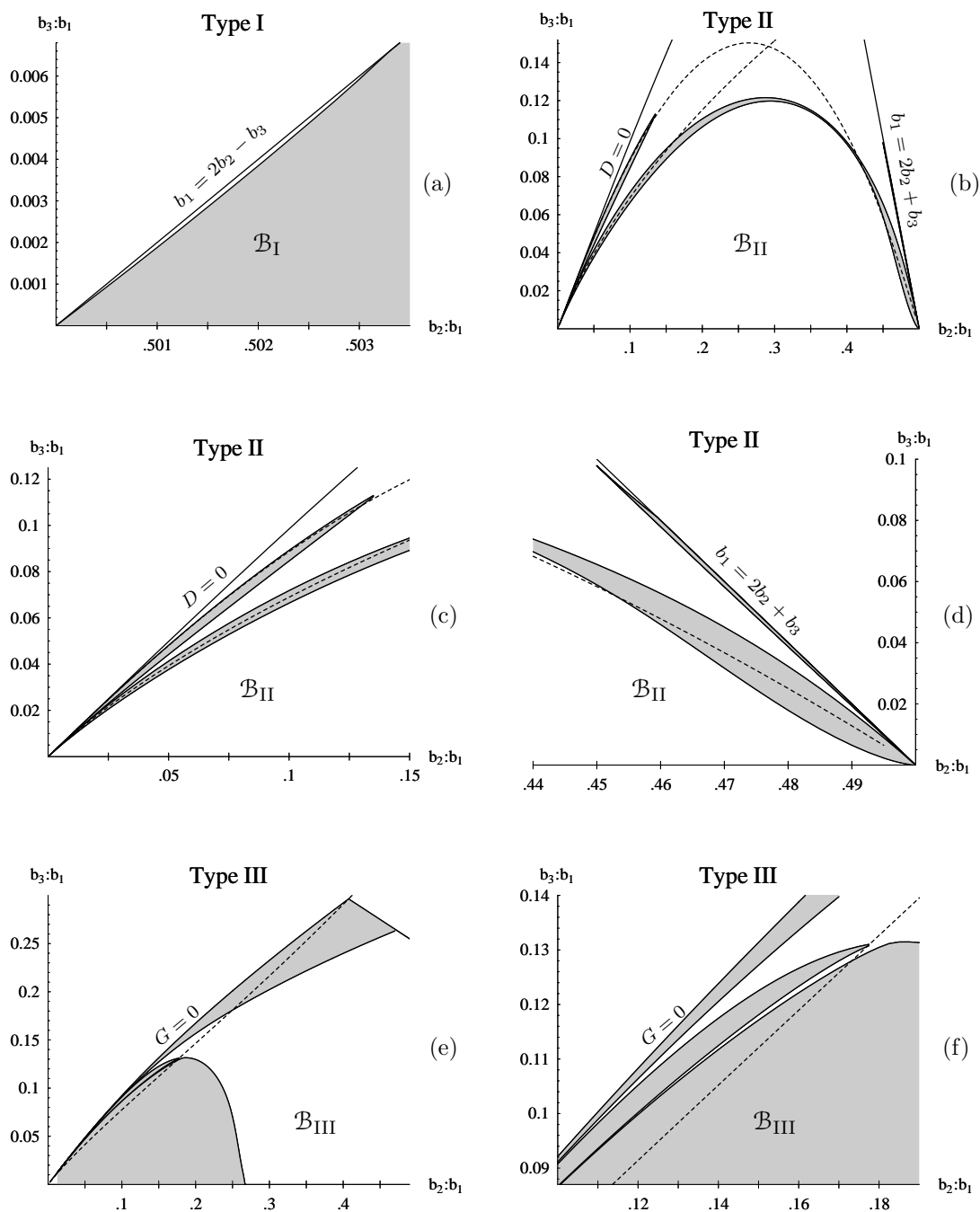}
{
\ins{69truemm}{162truemm}{\small (a)}
\ins{147truemm}{162truemm}{\small (b)}
\ins{44truemm}{153truemm}{$\mathrm{\cBI}$} 
\ins{115truemm}{153truemm}{$\cBII$} 
\ins{32truemm}{170truemm}{\rotatebox{43}{\scriptsize$b_1=2b_2-b_3$}}
\ins{95truemm}{169truemm}{\rotatebox{66}{\scriptsize$\D=0$}}
\ins{138truemm}{171truemm}{\rotatebox{-80}{\scriptsize$b_1=2b_2+b_3$}}
\ins{70truemm}{98truemm}{\small (c)}
\ins{146truemm}{98truemm}{\small (d)}
\ins{33truemm}{104truemm}{\rotatebox{40}{\scriptsize $\D=0$}}
\ins{112truemm}{107truemm}{\rotatebox{-45}{\scriptsize$b_1=2b_2+b_3$}}
\ins{42truemm}{87truemm}{$\cBII$} 
\ins{100truemm}{87truemm}{$\cBII$} 
\ins{70truemm}{34truemm}{\small (e)}
\ins{146truemm}{34truemm}{\small (f)}
\ins{29truemm}{39truemm}{\rotatebox{40}{\scriptsize $G=0$}}
\ins{101truemm}{39truemm}{\rotatebox{40}{\scriptsize $G=0$}}
\ins{50truemm}{23truemm}{$\cBIII$} 
\ins{118truemm}{23truemm}{$\cBIII$} 
}
\caption{\small Details of the regions of ellipticity (shaded) 
of the ellipsoids of types I, II and III. All curves are 
as in figure~\ref{Ell-1}.
\label{Ell-2}
}
\end{figure}

$\bullet$ The elliptic \Sdue ellipsoids of unknown Lyapunov
stability lie between the ``irrotational'' curve and
the border of region $\cBIII$. This is the dark shaded region in
figure~\ref{Ell-1}.b.

$\bullet$ Our computations confirm the existence of two
ellipticity regions for the ellipsoids of type~I, as shown in
figure~3 of \bib{chandrasekhar66}, but disclose a finer
structure of the lower region than that found by Chandrasekhar:
the upper border of this region appears to transversally
intersect the border $b_1=2b_2-b_3$ of the existence region at
about $b_2/b_1=.515$, but there is a very narrow crescent of
nonelliptic ellipsoids along this line, approximately in the range
$.5<b_2/b_1<.503$, which is shown in figure~\ref{Ell-2}.a. We could
not determine with certainty if this crescent actually touches
the line $b_1=2b_2-b_3$.

$\bullet$ The region of elliptic ellipsoids of type II consists
of three narrow fringes, shown enlarged in
figure~\ref{Ell-2}.b,c,d. None of these regions were detected by
Chandrasekhar, who stated that no ellipsoid of type II is
elliptic (\bib{chandrasekhar66}, pg. 169). Note that two of these
fringes belong (for the most part) to the region of positive 
relative pressure considered by Chandrasekhar. 

$\bullet$ The region of ellipticity of the ellipsoids of type III
has a three--lobed structure, the three lobes being separated by
two very thin gaps; see figure~\ref{Ell-2}.e, f. Chandrasekhar
detected only the upper of these three lobes. To within the
precision of our numerical computations, the two gaps close at
approximately $b_2/b_1=.01$, but it is possible that, in fact,
they persist to the origin.

In an attempt to understand the discrepancies between our results
and Chandrasekhar's, we tested our calculations using his
linearization of the equations of motion, finding complete
agreement. It is therefore possible that Chandrasekhar's
conclusions were simply based on too few sample cases to obtain
a detailed global picture of the structure of the regions of
ellipticity. A more detailed comparison with Chandrasekhar's
result is given in Appendix~B.

We conclude this analysis by observing that the sources of the
richness and intricacy of the structure of the ellipticity
regions resulting from these numerical computations are not yet
understood. Of course, we cannot be sure that our numerical
computations, even though accurate, ultimately detected the
finest details of this structure, particularly in view of the
nontrivial structure of the ellipticity regions for ellipsoids
approaching a degenerate (flat or rodlike) configuration. Many
features remain to be understood. For instance,
figures~\ref{Ell-1}.a and~\ref{Ell-2}.f suggest that there are
(nontrivial) relations between the ellipticity regions of the
ellipsoids I (resp. III) and the boundary of the existence
region of the ellipsoids \Sdue (resp. \Strep), which could be
the object of a bifurcation analysis. Additional analytical
insight into the problem is certainly needed.

\acapo
{\em Remarks: } (i) In his articles, Chandrasekhar refers to
ellipticity as ``stability'', while Riemann designated by this
name what is now called Lyapunov stability. Because of this, and
the previously mentioned fact that Riemann regarded as unstable
the ellipsoids which are saddle points of the Lyapunov function,
Riemann and Chandrasekhar found different ``stability'' and
``instability'' regions. Chandrasekhar interpreted this
difference as being due to some errors in Riemann's work (see
page 875 in \bib{chandrasekhar66} and pages 185--7 in
\bib{chandrasekhar69}). Riemann's results are fully confirmed by
our numerical computations (with the exception of the
characterization of the saddle points as unstable, of course)
and are, in fact, compatible with Chandrasekhar's results once
the distinction between Lyapunov stability and ellipticity is
taken into account.

(ii) In the quoted references, Chandrasekhar refers to
\cite{lebovitz66} for an explanation of ``the origin of
his [Riemann's] errors''. In that article, to the best of our
understanding, the source of these supposed errors is ascribed
to the fact that Riemann treated the Dirichlet problem as a
Lagrangian system (which is correct, as we have noted in Section
2), while according to \cite{lebovitz66} the system is
Lagrangian only in the case of the ellipsoids of type S; remark
(vi) in page 185 of \bib{chandrasekhar69} seems to indicate that
Chandrasekhar shared this opinion.

\acapon
{\bf B. Numerical Computations. } We now briefly describe
the numerical procedures used in the verification of
Numerical Conclusions 1 and in the Nekhoroshev analysis.
We also provide the key components of the proof of statement
(ii) of Proposition 5. We use canonical coordinates
obtained by introducing a pair of local ``Poincar\'e
coordinates'' $(q,p)$ on each sphere in the reduced phase
space; these coordinates coincide with the Poincar\'e elements
used in the Kepler problem. The use of coordinates is not
strictly necessary for the ellipticity analysis, but is very
useful for the construction of the normal forms used in the
analysis of the Nekhoroshev stability.

\def\Rad{\r}
\def\tm{\tilde m}
\def\tM{\tilde M}

\Lem{coordinates} For any $\Rad>0$, the map 
\beq{SPC}
   (q,p) \;\mapsto\; \tm_\Rad(q,p) \;:=\; 
     \left( p\,\sqrt{\Rad-{q^2+p^2\over4}} \,,\; 
          -\,q\,\sqrt{\Rad-{q^2+p^2\over4}} \,,\;
          \Rad- {q^2+p^2\over2} \right)
\eeq
is a diffeomorphism from the disk
$\{(q,p)\in\rdue:\,q^2+p^2<2\Rad\}$ onto the half--sphere
$\big\{m\in\rtre:\, \|m\|=\Rad,\, m_3>0\big\}$. It is symplectic
in the sense that $(\tilde m_\rho)_*\big(dp\wedge
dq\big)(m\per\o,m\per\o')= -m\cdot\o\per\o'$ for all $m$ in the
half sphere and all $\o,\o'\in\rtre$.
\eLem 

\acapon
{\bf Proof. } The map \for{SPC} is smooth and onto; its inverse
$m\mapsto \sqrt{2/(\Rad+m_3)}(-m_2,m_1)$ is also smooth. If
$v=(v_1,v_2,v_3)$ is tangent to the sphere at the point
$m=\tm_\Rad(q,p)$, then $v=v_q \der{\tm}{q} +v_p \der{\tm}{p}$
for some constants $v_q$ and $v_p$; specifically, if we set
$k_1=2\Rad-q^2-p^2$ and $k_2=k_1+2\Rad$, then $k_1 \sqrt{k_2}
v_q = qpv_1-(k_2-p^2)v_2$ and $k_1 \sqrt{k_2} v_p=(k_2-q^2)v_1-
qpv_2$. Thus $(\tilde m_\rho)_*\big(dp\wedge dq\big)(m\per\o,m\per\o')
  = (m\per \o)_p \,(m\per \o')_q - (m\per \o)_q \, (m\per 
\o')_p = - m\cdot(\o \per \o')$ for any $\o, \o' \in \rtre$.  \finelemma

\acapo
When studying a given equilibrium $(b^*,0,M^*)$, with
$M^*=(m^*_l,m^*_r)$, we shall translate the origin of the
coordinates $b$ to the equilibrium value $b^*$; for simplicity,
we will retain the notation $b$ for the translated coordinates.
Moreover, we shall use `rotated' Poincar\'e coordinates
$(q_1,p_1)$ and $(q_2,p_2)$ centered at $m^*_l$ and $m^*_r$ on
the spheres; specifically, for given $m^*\in\rtre$ we define
$$
    (q,p) \mapsto \tm(q,p;m^*) := R_{m^*}\; \tm_{\|m^*\|}(q,p)
$$
where $R_{m^*}$ is an orthogonal matrix such that
$R_{m^*}e_3={m^*\over\|m^*\|}$. Since rotations are symplectic,
the coordinates $\tilde m(\,\cdot\,;m^*)$ are symplectic. These
coordinates are not uniquely defined: right--multiplication of
$R_{m^*}$ by a rotation about $e_3$ produces a new system of rotated
Poincar\'e coordinates centered at $m^*$. To compress the
notation, we write $(q,p)$ for $(q_1,q_2,p_1,p_2)$
and\footnote{In what follows, to keep the notation simple, we
will write all vectors as row vectors, even though they must be
regarded as column vectors in all matrix formulas.}
$$
     M(q,p;M^*)\ug \big(\tm(q_1,p_1;m^*_l),
                        \tm(q_2,p_2;m^*_r)\big) \,.
$$
Thus, the reduced Hamiltonian \for{Hmu} takes the form
$$
   \cH(b,c,q,p;b^*,M^*) 
   \ug {1\over2}\,c\cdot \K(b^*+b)c \piu
   {1\over2}\,M(q,p;M^*) \cdot \J(b^*+b) M(q,p;M^*) \piu 
   \cV(b^*+b) \,.
$$

We now study the block structure of the Hessian matrix
$\cH''(0;b^*,M^*)$. To obtain a Hessian with as few nontrivial
blocks as possible, we construct the coordinates using matrices
such that $R_{m_l^*}e_1 = R_{m_r^*}e_1$ is parallel to a
principal axis orthogonal to $\mbox{span}[m_l, m_r]$. Using this
choice of coordinates and denoting by $\cH''_{bb}$ the
two--by--two block $\dder{\cH}{b_i}{b_j}$, etc, we have

\Lem{Hessian} For all equilibria, the blocks $\cH''_{cb}$,
$\cH''_{cq}$, $\cH''_{cp}$, $\cH''_{bp}$, and $\cH''_{qp}$ equal
zero and $\cH''_{cc}=\K$.
Moreover, $\cH''_{bq}=0$ for the equilibria of type S. 
\eLem

\acapon
{\em Proof. } 
The equalities $\cH''_{cb} = \cH''_{cq} = \cH''_{cp} = 0$ and
$\cH''_{cc}=\K$ are obvious. Note that the `standard' Poincar\'e 
parameterization \for{SPC} of the sphere satisfies
\beq{der-poi-1}
  \der \tm q(0,0;\Rad)  \ug -\sqrt\Rad\,e_2 \,, \qquad
  \der \tm p(0,0;\Rad)  \ug  \sqrt\Rad \,e_1 \,,\qquad
  \dder\tm qp(0,0;\Rad)\ug 0 \,.
\eeq
Using the latter equality and denoting by $\der{\tilde
M}{q_j}$ (resp. $\der{\J}{q_j}$) the vector (resp. matrix) with 
entries given by the derivatives with respect to $q_j$
of the entries of $\tilde M$ (resp. $\J$), etc., we have 
\begin{eqnarray}
   \dder \cH{b_i}{p_j}(0;b^*,M^*) \ugarr 
       R\, \der{\tilde M}{p_j}(0,0) 
       \cdot \der{\J}{b_i}(b^*) \,M^*       \nonumber\\  
   \dder \cH{q_i}{p_j}(0;b^*,M^*) \ugarr
       R\, \der{\tilde M}{q_i}(0,0) 
       \cdot \J(b^*) R\, \der{\tilde M}{p_j}(0,0) \label{owl} \\
   \dder \cH{b_i}{q_j}(0;b^*,M^*) \ugarr 
       R\, \der{\tilde M}{q_j}(0,0) 
       \cdot \der{\J}{b_i}(b^*) \,M^*  \,,      \nonumber
\end{eqnarray}
where $R$ is the six--by--six block diagonal matrix with 
diagonal blocks $R_{m_l^*}$ and $R_{m_r^*}$.

Due to the block structure of $\J$, the two `components' of
$\der{\J}{b_i}(b^*)M^*$ lie in $\mbox{span}[m_l, m_r]$. (Here 
we refer to $m_l$ and $m_r$ as the `components' of $(m_l,m_r)$.) 
On the other hand, \for{der-poi-1} and the condition that $R_{m_l^*} e_1
= R_{m_r^*} e_1$ be orthogonal to $\mbox{span}[m_l, m_r]$ imply 
that the components of $R\der{\tilde M}{p_j}(0,0)$ are orthogonal
to $\mbox{span}[m_l, m_r]$ for $j = 1, 2$; hence \for{owl} implies  
that $\dder \cH{b_i}{p_j}(0;b^*,M^*)=0$. To see that 
$\dder \cH{q_i}{p_j}(0;b^*,M^*)$ equals zero, note
that $R_{m_l^*} e_1 = R_{m_r^*} e_1$ is, by construction, an
eigenvector of $J_1(b^*)$ and $J_2(b^*)$; hence \for{der-poi-1} 
implies that the components of $R^T\J(b^*) R\, \der{\tilde M}{p_j}(0,0)$ 
are orthogonal to the components of $\der{\tilde M}{q_j}(0,0)$.

In the case of the S ellipsoids, the two components of $M^*$ are
both parallel to the same coordinate axis. Due to the block
structure of $\J$, the components of $\der{\J}{b_i}(b^*)M^*$ are
parallel to the components of $M^*$; thus, both components of
$R^T \der{\J}{b_i}M^*$ are parallel to $e_3$. By \for{der-poi-1}
and \for{owl}, this implies $\cH''_{bq}=0$. 
\finelemma 

\acapo
Let us now consider the ellipsoids of type \Sdue. With a suitable
choice of the coordinates, the Hessian has the diagonal block
structure 
$$
   \cH''(0;b^*,M^*) \ug 
   \diag \big[\cH''_{bb}, \K, \cH''_{qq}, \cH''_{pp} \big]\,.
$$
We limit ourselves to stating without proof the following facts,
which can be proven using elementary, though sometimes laborious, 
computations:
\bList
\item[$\bullet$] The blocks $\K$ and $\cH''_{qq}$ are positive 
definite for all the \Sduep--ellipsoids.
\item[$\bullet$] The block $\cH''_{pp}$ is positive definite for all
the counterparallel \Sduep--ellipsoids and is indefinite for all 
the coparallel \Sduep--ellipsoid. 
\end{list}
\noindent
The block $\cH''_{bb}$ has a relatively simple expression, but
the analytical study of its definiteness is quite laborious. 
We verified numerically that it is positive definite at all 
ellipsoids of type \Sdue. (Chandrasekhar's study of this
block in \bib{chandrasekhar65} is also numerical.)
The linearization at the equilibrium of the Hamiltonian vector
field of $\cH$ is $X(b^*,M^*) \big(b,c,q,p\big)$, where the
matrix $X(b^*,M^*)$ is the product of the symplectic matrix and
$\cH''(0;b^*,M^*)$. Thus, at a type \Sdue ellipsoid, the matrix
$X(b^*,M^*)$ is block diagonal, with the two $4\per4$ blocks
$$
   X_{bc}
   \ug \left(\matrix{  0                 &\K   \cr
                      -\cH''_{bb} &0  }\right)
   \quad\mathrm{and}\quad
   X_{qp}
   \ug \left(\matrix{  0          &\cH''_{pp}\cr
                      -\cH''_{qq}     &0          }\right) \,.
$$ 
Since the eigenvalues of each of these two matrices are purely 
imaginary if its two nonzero blocks are positive semi--definite,
we arrive to statement (i) of Numerical Conclusions~1. This 
conclusion is not rigorous because it is based on the numerical
verification of the definiteness of $\cH''_{bb}$. Statement (ii)
of Proposition~\ref{LyapStab} follows from the definiteness of
$\cH''_{qq}$ and the indefiniteness of $\cH''_{pp}$.

\acapo
{\em Remark: } In \bib{riemann860}, Riemann suggests an analytic
argument that would show that for all ellipsoids of type S
satisfying $m_l \neq m_r$, i.e. $\MS_-(b_1, b_2, b_3) > 0$, the
function $b\mapsto\cH(b,0,0,0;b^*,M^*)$ has a minimum at $b=0$.
Since the block $\cH''_{bb}$ is positive semi--definite at such a
minimum, this result could be used to construct an analytic
proof of statement (i) of Numerical Conclusions~1. However, as
Riemann explicitly states, he does not actually carry out the
minimality argument in \bib{riemann860}.

\acapo
We now consider the ellipsoids of types I, II and III. The
Hessian has the form
$$
   \cH''(0;b^*,M^*) \ug 
   \left(\matrix{ 
    \cH''_{bb}  &0      &\cH''_{bq}  &0         \cr
        0      &\K        &0         &0         \cr
    \cH''_{qb}  &0      &\cH''_{qq}  &0         \cr
        0       &0      &0           &\cH''_{pp}  } 
   \right) \,.
$$
The eigenvalue problem of this matrix factorizes into the
eigenvalue problems for the three matrices
$$
    \K \,,\qquad \cH''_{pp} \,,\qquad
     \left(\matrix{ 
          \cH''_{bb} &\cH''_{bq} \cr
          \cH''_{qb} &\cH''_{qq} }\right) \,.
$$
One can prove that the two--by--two matrix $\cH''_{pp}$ is 
indefinite for any ellipsoid of type I, II or III. (This
implies that these ellipsoids are saddle points of the reduced
Hamiltonian, as stated by Riemann.)
The eigenvalue problem for the linearization of the Hamiltonian
vector field determined by this Hessian does not factorize.
Therefore the analysis of the ellipticity is particularly
difficult. We know of no rigorous ellipticity results for
these ellipsoids. 

\acapon

\newsection{Nekhoroshev--stability of the Riemann--ellipsoids: 
Generalities}

We base our study of the Nekhoroshev--stability of the Riemann
ellipsoids on some recent results in Nekhoroshev theory, which
we describe in this section with reference to the problem at
hand; complete proofs can be found in \bib{benfasguz98} (see also
\bib{fasguzben98}, \bib{niederman98}, \bib{guzfasben98}). As is
common in Hamiltonian perturbation theory, the procedure consists of 
the construction of the Birkhoff normal form of order four for the 
reduced Hamiltonian and testing of the normal form for certain
properties that we now specify.

\acapo

\acapon
{\bf A. Birkhoff Normal Form. } Consider an elliptic equilibrium
$\big(b^*,M^*\big)$ with $M^*=(m_l^*,m_r^*)$,
$m^*_l\not=0\not=m^*_r$. Use the canonical coordinates
$\xi=(b,c,q,p)$ defined as in Section~4.B. Let
$\cH_j(\xi;b^*,M^*)$ denote the $j$--th term in the Taylor
series at $\xi=0$ for the reduced Hamiltonian $\cH$, so that
$\cH_j$ is a homogeneous polynomial of order $j$ in $\xi$ and
$$
  \cH(\xi;b^*,M^*) = \cH_2(\xi;b^*,M^*) +\cH_3(\xi;b^*,M^*)  
                  + \cH_4(\xi;b^*,M^*) \piu \ldots.
$$
Let $X(b^*,M^*)\xi$ be the Hamiltonian vector field of $\cH_2$.
Since the equilibrium is elliptic, the matrix $X(b^*,M^*)$ has
purely imaginary eigenvalues
$$
     \pm i \, \o_j(b^*,M^*) \,,\qquad j=1,\ldots,4 \,,
$$
with the convention $\o_j\ge 0$.

The first step in the construction of the Birkhoff normal
form for the Hamiltonian $\cH$ is the ``symplectic
diagonalization'' of the quadratic term $\cH_2$. It is known
that if the eigenvalues of $X(b^*,M^*)$ are all distinct, then
there is a linear canonical change of coordinates
$$
    \xi \mapsto \Xi=(B,C,Q,P)\ug T(b^*,M^*)^{-1} \xi 
$$
such that $\widehat\cH_2(\Xi;b^*,M^*):=\cH_2(T(b^*,M^*)\Xi;b^*,M^*)$ 
has the form
$$
    \widehat\cH_2\ug  
    s_1 \o_1{{B_1^2+C_1^2}\over 2} \piu 
    s_2 \o_1{{B_2^2+C_2^2}\over 2} \piu 
    s_3 \o_3{{Q_1^2+P_1^2}\over 2} \piu 
    s_4 \o_4{{Q_2^2+P_2^2}\over 2} \,, 
$$
where $s_j=\pm1$. The matrix $T$ and the numbers $s_j$ 
are constructed as follows: For $j = 1,..,4$, let $x_j^\pm = x_j' \pm
i x_j'' \in \complessi^8$, with $x_j',\,x_j'' \in \reali^8$,
denote any eigenvectors of the matrix $X(b^*,M^*)$ associated to
the eigenvalues $\pm i\o_j$. One can show that 
$$
   \Gamma_j(b^*,M^*) \,:=\; x_j'\cdot J_8 x_j'' \;\not=\; 0\,, 
   \qquad j=1,..,4 \,,
$$
where $J_8 = \diag\big[J_4,J_4\big]$; here $J_4$ denotes the
standard four--by--four symplectic matrix $\left(\matrix{\mathbf{0} 
&\unity \cr -\unity &\mathbf{0} }\right)$. We set
$(u_j,v_j)=(x_j',x_j'')$ if $\Gamma_j>0$ and
$(u_j,v_j)=(x_j'',x_j')$ if $\Gamma_j<0$, and define 
$$
    T \ug \Big( 
      {u_1\over\g_1},{u_2\over\g_2},{v_1\over\g_1},{v_2\over\g_2},
      {u_3\over\g_3},{u_4\over\g_4},{v_3\over\g_3},{v_4\over\g_4}\Big) 
    \,,\qquad       
    s_j \ug \mathrm{Sign}(\Gamma_j) \,,
$$
where $\g_j=\sqrt{|\Gamma^j|}$; the notation means that the eight
vectors are the rows of the matrix. Note that the condition that
all the eigenvalues of $X$ be distinct (and hence nonzero) can
be regarded as the condition that the ``frequency vector''
$$
      \O \ug (s_1\o_1,\ldots,s_4\o_4)
$$
has no resonances of order one or two:
\beq{res-2}
    \O\cdot\nu \;\not=0\; \quad\mathrm{for\ all\ } 
    \nu\in\interi^2 \quad\mathrm{such\ that\ } 
    |\nu|=1,2,
\eeq
where, for integer vectors, $|\nu|=|\nu_1|+\cdots+|\nu_4|$.

\acapo
{\it Remark: } In the case of the ellipsoids of type S, the
eigenvalue problem for the matrix~$X(b^*,M^*)$ factorizes into
that for two $4\per4$ blocks. Since the construction of the
matrix $T$ also factorizes, it is necessary to check the
nonresonance condition \for{res-2} for all resonances of order
one, but only for the four resonances $(1,\pm1,0,0)$,
$(0,0,1,\pm1)$ of order two.

\acapo
For the construction of the Birkhoff normal forms we 
use the complex coordinates $U = 
(W_1,W_2,W_3,W_4,Z_1,Z_2,Z_3,Z_4)$ defined by
$$
   W_j = - {iB_j+C_j\over\sqrt2} \,,\quad 
   Z_j = {B_j+iC_j\over\sqrt2} \,,\quad 
   W_{j+2} = - {iQ_j+P_j\over\sqrt2} \,,\quad 
   Z_{j+2} = {Q_j+iP_j\over\sqrt2} \quad (j=1,2)
$$
(note that $W$ are the coordinates, $Z$ the momenta). Since 
$\Xi = \Sigma U$, with
$$
  \Sigma=\left(\matrix{\Sigma_4 & 0 \cr 
                       0 & \Sigma_4 \cr} \right)
  \qquad\mbox{and}\qquad 
  \Sigma_4={1\over\sqrt2}\left(\matrix{
            i& 0&1 &0 \cr
            0& i& 0&1\cr
           -1& 0& -i& 0\cr
            0&-1& 0&-i \cr} \right) \,,
$$  
the Hamiltonian in terms of the complex coordinates $U$ is
$H(U;b^*,M^*) :=\cH(T \Sigma U ;b^*,M^*)$ and its
quadratic part is
$$
    H_2(U;b^*,M^*) \ug \sum_{j=1}^4 i\O_j(b^*,M^*)Z_jW_j \,.
$$

We now describe the construction of the Birkhoff normal form of
order four of the Hamiltonian~$H$. Given a function $f$ of
$W=(W_1,W_2,W_3,W_4)$ and $Z=(Z_1,Z_2,Z_3,Z_4)$, we denote its 
Taylor series by $f(W,Z)=\sum_{j,k\in \naturali^4}f_{jk}W^jZ^k$,
where $W^j = W_1^{j_1}W_2^{j_2}W_3^{j_3}W_4^{j_4}$ etc.
For any integer vector $\nu\in\interi^4$, we define the
$\nu$--th {\em harmonic} of $f$ by
$$
   \fou f\nu(W,Z) \ug 
   \sum_{j,k\in \naturali^4\atop j-k=\nu}f_{jk}W^jZ^k \,.
$$
The {\em average} of $f$ is its harmonic $\fou f0$. The
{\em spectrum} of $f$ is $\Sp f=\{\nu\in\interi^4: \fou
f\nu\not=0\}$. 

We construct the Birkhoff normal form by means of the so--called
Lie method. Let $\P^\chi$ denote the time--one--map of the flow
of the Hamiltonian vector field of a function $\chi$; the Lie
transform generated by $\chi$ of a function $f$ is
$f\circ\Phi^\chi$. If $f$ and $\chi$ are analytic in an open
neighbourhood of a point (the equilibrium point, in our case),
then $f\circ\Phi^\chi$ is also analytic in some sufficiently
small, but nonempty, open neighbourhood of such a point and
there one has $f\circ\Phi^\chi=\sum_{j=0}^\infty {1\over
j!}\,L_\chi^j f$. (We adopt the sign convention $L_f g = \{f,g\}
=\der fZ\der gW -\der gZ\der fW$ for the Poisson brackets.) 

The generating function $\chi$ of the Lie transform is
constructed by solving the so--called homological equation,
which has the form $\{H_2,\chi\}=f-\fou{f}0$ for some function
$f$. The solution of this equation is formally given by
$\chi=\cS_\O(f)$, where
$$
  \cS_\O(f) \ug  \sum_{\nu\in\Sp f \setminus\{0\} } 
               {\fou{f}\nu\over i\O \cdot \nu} \,.
$$ Obviously, this is well defined if $\O$ does not resonate with
any $\nu\in \Sp f$, and if the series converges. In our case $f$
will always be a polynomial, either $f=H_3$ or $f=H_4'$ (see
below), so no convergence problems exist.

The fourth--order Birkhoff normal form of $H=H_2+H_3+H_4+\ldots$,
if it exists, is constructed using two Lie transforms, which
average the terms of degree three and four of $H$, respectively.
The first Lie transform is generated by the solution 
\beq{chi1}
  \chi_1 \ug \cS_\O(H_3)
\eeq
of the homological equation $\{H_2,\chi_1\}=H_3$ and is well
defined if
\beq{res-3}
    \O(b^*,M^*)\cdot\nu\not=0 \qquad \forall\nu\in\Sp{H_3} \,.
\eeq
Since $\fou{H_3}0=0$ and $L_{\chi_1}^2H_2=-L_{\chi_1}H_3$ one 
sees that this Lie transform conjugates $H=H_2+H_3+H_4+....$ 
to 
$$
    H' \ug H_2 \piu H'_4 \piu \ldots \,,\qquad \mbox{where} \qquad
    H'_4 \ug {1\over2}L_{\chi_1} H_3 + H_4 \,.
$$
The second Lie transform is generated by the solution
$\chi_2 = \cS_\O(H'_4)$ of $\{H_2,\chi_2\}=H'_4-\fou{H'_4}0$. It is
well--defined if 
\beq{res-4}
    \O(b^*,M^*)\cdot\nu\not=0 \qquad 
    \forall\nu\in\Sp{H'_4}\setminus\{0\} 
\eeq
and leads to the fourth--order Birkhoff normal form
\beq{BNF}
     H'' \ug H_2+H''_4+\ldots 
     \qquad \mbox{with} \qquad
     H''_4=\fou{\half L_{\chi_1} H_3 + H_4}{0} \,.
\eeq
For our purposes, the only important points are controlling the
existence of the fourth--order Birkhoff normal form and testing
that it possesses the convexity properties discussed above. 
Neither of these tasks require the actual computation of the
generator $\chi_2$: one needs only test \for{res-4} and compute 
$H''_4$. 

\acapon

\acapon 
{\bf B. Nekhoroshev stability of elliptic equilibria. } Let us
introduce the action functions $I=(I_1,\ldots,I_4)$, 
$I_j=iW_jZ_j$, and regard $H''_4$ as a quadratic form on
$\reali^4$ by writing 
$$
    H''_4(I) \ug {1\over2}I\cdot A(b^*,M^*)I 
$$
where $A(b^*,M^*)$ is a $4\per4$ symmetric matrix. The Birkhoff
normal form of order four \for{BNF} is said to be
\bList
\item[$\bullet$] {\it convex} if the quadratic form $H''_4$ is definite, 
i.e. the eigenvalues of $A$ are either all positive or all negative.
\item[$\bullet$] {\it quasi--convex} if the restriction of $H''_4$ 
to the subspace orthogonal to $\O$ is definite, i.e. 
$$
    \O\cdot I =0 \quad\mathrm{and}\quad H''_4(I)=0 
    \quad\Longrightarrow\quad  I=0 \,.
$$
\item[$\bullet$] {\em directionally quasi--convex} if the
restriction of the quadratic form $H''_4$ to the plane
orthogonal to $\O$ is nonvanishing in the ``first 
16--ant'', i.e. 
$$
   \O\cdot I \ug 0 \,,\quad H''_4(I)\ug 0
   \quad\mathrm{and}\quad I_1,\ldots,I_4\ge0
   \qquad\Longrightarrow\qquad  I=0 \,.
$$
\end{list}
Clearly, each notion generalizes the previous one.

If the Hamiltonian $H(W,Z)$ and its Birkhoff normal form of order
four are analytic, then any of these conditions implies the
Nekhoroshev--stability of the equilibrium, in the precise sense
that {\em for any small $\e>0$ one has
$$
   \|I(0)\|\le\e \quad\Longrightarrow\quad \|I(t)\|\le\e^{\a}
   \quad \mbox{\rm for} \quad |t|\le \exp \e^{-\beta}
$$
for some positive constants $\a$ and $\beta$}. There are different 
possible values for the constants $\a$ and $\beta$, which we report 
here with reference to the case under consideration of a system 
with four degrees of freedom:
\bList
\item[--] If the Birkhoff normal form is directionally
quasi--convex, then one can take $\a=\b={1\over4}$ and it is
also possible to prove a stricter confinement of motions on
correspondingly shorter times, namely $\a={1+k\over4+k}$,
$\b={1\over4+k}$ for any $k>0$ (see \bib{fasguzben98}).
\item[--] If the Birkhoff normal form is quasi--convex
(directional quasi--convexity does not seem to be sufficient),
then it is possible to prove the above estimates with $\a=1$ and
$\b={1\over16}$. (See \bib{guzfasben98}.)
\item[--] If the Birkhoff normal form of order four is quasi--convex 
and it is possible to construct the Birkhoff normal form of order 
$s>4$, then one also has $\a=1$, $\b={s-3\over16}$ (see 
\bib{guzfasben98}). 
\end{list}

\acapo 
{\em Remark: } The original motivation for the quoted articles
\bib{fasguzben98}, \bib{guzfasben98} on the Nekhoroshev
stability of elliptic equilibria was precisely the present study
of the Riemann ellipsoids. Until a few years ago, all Nekhoroshev
stability results assumed that the frequency vector would
satisfy a strong nonresonance condition, typically a diophantine
one. This would be insufficient for the study of a problem such
as the present one, in which the frequencies of the equilibrium,
and therefore its nonresonance properties, depend on
continuously varying parameters. The notion of directional
quasi--convexity was introduced in \bib{benfasguz98} in
connection with the study of the Nekhoroshev--stability of the
triangular points of the restricted three--body problem. As will
be seen below, this notion also plays a central role in the
study of the Riemann ellipsoids, since the majority of the
ellipsoids are not quasi--convex. (Nekhoroshev estimates also
hold if a sufficiently high order Birkhoff normal form exists
and satisfies some ``steepness'' condition; directional
quasi--convexity allows us to avoid these generalizations, the
hypotheses of which would be rather difficult to verify and
which would lead to significantly worse values for the constants
$\a$ and $\beta$.)

\acapon

\acapon 
{\bf C. KAM theory. } If the Birkhoff normal form of order four
is {\it nondegenerate}, in the sense that $\det
A(b^*,M^*)\not=0$, then KAM theory applies, ensuring that in any
sufficiently small neighbourhood of the equilibrium, the
majority of the initial data gives rise to motions which are
quasi--periodic with four frequencies. (See e.g. \bib{arnold88}
\bib{meyer+hall91}.)

\acapo

\newsection{Nekhoroshev--stability of the Riemann ellipsoids: 
Numerical results} 

{\bf A. Results. } We have numerically constructed the Birkhoff
normal forms for a (quite large) number of sample Riemann
ellipsoids. Specifically, within the existence region of each
type of ellipsoid, we considered a mesh in the plane
$\big({b_2\over b_1}, {b_3\over b_1}\big)$ determined by
vertical lines uniformly spaced at a distance $.0025$; the
number of mesh points on each vertical lines typically varied
between twenty and one hundred, depending on the length of the
line (however, in some cases we used as many as five hundred
points). This analysis leads to the following

\acapon
{\bf Numerical Conclusions 2. } {\em The Birkhoff normal form of
order four exists and is directionally quasi--convex and
nondegenerate for all the Riemann ellipsoids except for those
lying on a finite number of curves in the set $\cB$ corresponding 
to resonances of order up to four.}
 
\acapo
We regard the fact that all of the computed nonresonant Birkhoff 
normal forms are directionally quasi--convex as indicating that all
nonresonant Riemann ellipsoids are Nekhoroshev stable. We
investigated the existence of low--order resonances in greater
detail, without restriction to the points of the considered mesh
(see subsection 6.B for some details), finding that there are a
finite number of curves in the plane $\big({b_2\over b_1},
{b_3\over b_1}\big)$ at which at least one of the nonresonance
conditions \for{res-2}, \for{res-3}, \for{res-4} is violated.
Specifically, we found 8 different resonances for the type
\Sduep, 52 for the type I, 33 for the type II, and 47 for the
type III. For example, figure~\ref{QC-Res}.a reports the
resonant curves for the ellipsoids of type~I.

We found that the Birkhoff normal form is quasi--convex only for
a few ellipsoids of each type; in turn, very few of these
ellipsoids are actually convex. For example,
figure~\ref{QC-Res}.b shows the quasi--convex ellipsoids in the
lower region of elliptic type I ellipsoids. The determination of
the quasi--convex ellipsoids has some interest because, as
discussed in the previous Section, such equilibria may have
stronger stability properties if they are nonresonant up to
sufficiently high order. (From this perspective, the finer
distinction between convexity and quasi--convexity does not
instead seem to be as significant.)

Nondegeneracy of the Birkhoff normal forms implies that the KAM
theorem applies, so that the majority of the motions near each
equilibrium in the reduced phase space are quasi--periodic with
four frequencies. ``Reconstruction'' of these motions then gives
quasi--periodic motions of the Dirichlet problem, with up to
(and in fact, typically) eight frequencies.

\begin{figure}
\vskip7.5truecm
\includegraphics{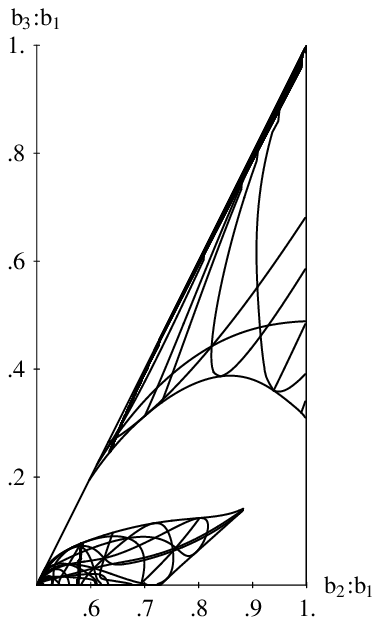}
\includegraphics{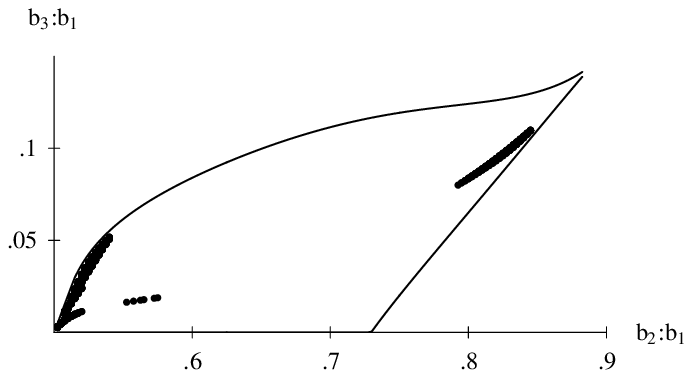}
\ins{6truemm}{36truemm}{\small (a)}
\ins{146truemm}{36truemm}{\small (b)}
\caption{\small (a) The ellipsoids of type I that are resonant
up to order four and (b) the quasi--convex Riemann Ellipsoids of
type I in the lower region of ellipticity.}
\label{QC-Res}
\end{figure}

\acapo

\acapon
{\bf B. Numerical Procedures. } We give now some information
about the numerical procedures adopted here. The fourth--order
Taylor expansion of the Hamiltonian in coordinates $(b,c,q,p)$,
and hence the Hamiltonian matrix $X(b^*,M^*)$, has been
constructed analytically for each type of ellipsoid. The
coefficients of these polynomials have been evaluated for any
ellipsoid in our mesh, and all other operations have then been
performed on polynomials with numeric coefficients. These
include:
\bList
\item[--] Computing the eigenvalues and the eigenvectors of
$X(b^*,M^*)$ and testing for the presence of resonances of order
one and two. Note that once the eigenvectors are known, the matrix
$S(b^*,M^*)$ determining the coordinates $\Xi$ and $U$ is
also known.
\item[--] Constructing the normal forms by means of
formulas \for{BNF} and \for{chi1}. This requires computing
Poisson brackets, i.e. derivatives, and averages of polynomials
with numeric coefficients; these are easily implemented operations.
\item[--] Checking the nonresonance conditions \for{res-3}
and \for{res-4}, which assure that the normal forms can be 
constructed. This involves determining the spectrum of 
a polynomial.
\item[--] Testing the normal forms, namely the numeric vector
$\O$ and matrix $A$, for convexity, quasi--convexity, or
directional quasi--convexity.
\end{list}
The chosen implementation of the latter two points requires
some explanation.

\acapo
{\em Spectra and Resonances. } We computed the spectra of $H_3$
and $H''_4$ at every point $(b^*,M^*)$ in the considered mesh
and took the unions $\Sigma_3$ and $\Sigma_4$ of all these sets.
(The union is taken because the numerically computed spectrum at
a given point may lack very small harmonics.) We identified as
resonances those $\nu\in \Sigma_3\cup \Sigma_4$ for which the scalar
product $\O(b^*,M^*)\cdot\nu$ either vanishes or changes sign on
our grid; of course, use of a finer grid could reveal additional
resonances. (The resonant curves of figure~\ref{QC-Res}.a were
constructed by computing the approximate zeros of
$\O(b^*,M^*)\cdot\nu$, testing for changes of sign on a finer
grid.)

\acapo
{\em Quasi and Directional Quasi--convexity. } If the matrix $A$
is not convex, then the tests for quasi--convexity and
directional quasi--convexity can be carried out as follows.
Let $R\in\SOTRE$ be such that $R\O=(1,0,0,0)$, so that the
Hessian of the restriction of the quadratic form
$H''_4$ to the subspace orthogonal to $\O$ is the lower
right hand $3\per3$ block $\tilde A$ of the matrix $RAR^T$.
Quasi--convexity of $H''_4$ is equivalent to convexity of the
matrix $\tilde A$. 

Under the additional hypothesis that the restriction of
$H''_4(I)$ to the subspace orthogonal to $\O$ is
nondegenerate,\footnote{This hypothesis is not crucial, but it
is satisfied by all the computed ellipsoids.} the test for
directional quasi--convexity can be performed as follows. The
quadratic form $H''_4$ is directionally quasi--convex iff any
``asymptotic vector'', namely any nonzero vector $I$ satisfying $I\cdot
\tilde AI=0$, points out of the first 16--ant, i.e. has at least
one negative and one positive entry. Since these vectors form a
cone (of dimension two, if $\tilde A$ is nondegenerate) we need
only test the asymptotic vectors in the intersection of the cone
and the unit sphere, namely two ellipses. In order to determine
these ellipses, note that in the absence of quasi--convexity,
the matrix $\tilde A$ has at least one eigenvalue of either
sign. We may assume that the eigenvalues of $\tilde A$ satisfy
$\a_1>0$ and $\a_2,\a_3<0$, replacing $\tilde A$ with $- \tilde
A$ if necessary. (This is possible because $\tilde A$ and $-
\tilde A$ have the same asymptotic vectors. Also, the
nondegeneracy hypothesis assures that all of the eigenvalues are
nonzero.) Let $\lcb \tilde v_1, \tilde v_2, \tilde v_3 \rcb$ be
an orthonormal eigenbasis associated to $\tilde A$, so that
$S =(\tilde v_1,\tilde v_2,\tilde v_3)\in\SOTRE$ satisfies 
$S^TAS = \diag(|\a_1|,-|\a_2|,-|\a_3|)$. In this basis, the cone 
of asymptotic vectors is given by
$$
   \Big\{ x=(x_1,x_2,x_3)\in\rtre: 
   x_1^2 \ug {|\a_2|\over\a_1} x_2^2 + {|\a_3|\over\a_1} x_3^2   
   \Big\} \,.
$$
The two ellipses of unit vectors on this cone consist of the vectors
$$
   \hat x^\pm(\theta) \ug \left(
   \pm \sqrt{ {|\a_2|\over \a_1+|\a_2|} \cos^2\theta +
              {|\a_3|\over \a_1+|\a_3|} \sin^2\theta } ,\,
   \sqrt{\a_1\over \a_1+|\a_2|} \cos\theta \,, 
   \sqrt{\a_1\over \a_1+|\a_3|} \sin\theta \right) \,,\quad
   \theta\in S^1 \,. 
$$
In the original basis in the subspace orthogonal to $\O$,
these vectors are 
$$
    \hat I^\pm(\theta) = (\hat I_1,\hat I_2,\hat I_3) 
    \ug S\hat x^\pm(\theta) \,.
$$
These are three--dimensional vectors belonging to the subspace
orthogonal to $(1,0,0,0) = R \O$, so the corresponding
four--dimensional vectors are $I^\pm(\theta)=R^T(0,\hat I_1,\hat
I_2,\hat I_3)$. The normal form is directionally quasi--convex
iff for every $\theta\in S^1$ the vector $I^+(\theta)$ has at
least one negative and one positive entry. (It is obviously
sufficient to consider only one of the two circles.) In the
computations, we have tested this condition by checking that, at
each of a number (typically $101$) of uniformly spaced values of
$\theta$ in the interval $0\le\theta\le2\pi$, there are two
eigenvalues whose product is negative.

\acapo

\newsection{Conclusions}

In this work, we have posed the problem of the Nekhoroshev
stability of those Riemann ellipsoids that are elliptic but of
unknown Lyapunov stability. We have provided numerical evidence
that all Riemann ellipsoids that are not resonant up to order
four are Nekhoroshev stable. We have developed a consistent and
rigorous Hamiltonian formulation of the problem on a
four--to--one covering manifold, which validates the procedures
(and clarifies the presence of a previously not understood
discrete symmetry in the problem). Within this formulation we
provided explicit, precise descriptions of the existence regions
and of the momenta of the Riemann ellipsoids of the various
types and we reviewed the existent results on the stability
and ellipticity of the Riemann ellipsoids, showing (numerically)
that some of the latter regions are significantly larger than was
previously thought.

At the conclusion of this work, it is our opinion that
several key features of the problem are not yet understood.
Some open problems include:
\bList
\item[--] Properties of axisymmetric Riemann ellipsoids other than the
Maclaurin spheroids. 
\item[--] Analytic determination of the regions of ellipticity. 
Understanding of their subtle structure. 
\item[--] Lyapunov stability or instability of those Riemann ellipsoids 
that are saddle points of the reduced Hamiltonian. 
\end{list}
\noindent
Of course, it would also be interesting to prove the
Nekhoroshev--stability rigorously, that is, without recourse to
numerical calculations, but this appears to be a formidable task
at present.

In conclusion, we would like to mention a point of historic
interest. Riemann's treatment of the Dirichlet problem is
exceptionally advanced and uses several geometric constructions,
including Poisson reduction, the amended potential, and the
trivialization of the tangent and cotangent bundles of Lie
groups, that are now considered to be central to modern
geometric mechanics.

Riemann's apparent lack of concern with proving the actual
instability of the saddle points of the reduced Hamiltonian
is probably typical of that era --- this error permeated the
mechanical literature long after Lyapunov clarified the subject
at the beginning of the twentieth century. At the dawn of a new
century, the Lyapunov (in)stability of equilibria that are both
elliptic and saddles of the `obvious' Lyapunov function is still
a formidable problem. Nekhoroshev stability provides a weaker,
but very practical, notion of stability, which has the ultimate
advantage of being compatible with Lyapunov instability.

\acapon

\newsection{Appendix A: On the proof of Proposition \ref{EqRel} }

The proof of Proposition~\ref{EqRel}, even though 
elementary, is too lengthy to be reported here. Hence we 
indicate the main steps in the form of Lemmas that we do not
prove. The relative equilibria are the critical points of the
reduced Hamiltonian $H_\eta:P_\eta\to\reali$ and hence are given
by the conditions
\beq{critical}
  \nabla_b H_\eta=0\,,\quad 
  \nabla_c H_\eta=0\,,\quad 
  m\per \nabla_m H_\eta =0\,.
\eeq
Since $\nabla_c H_\eta=0$ iff $c=0$, we can equivalently 
search for critical points $(b,m)$ of
$$
   \tilde H_\eta(b,m) \;:=\; 
   \half m\cdot \J(b)m \piu \cV(b) \,.
$$
Note that only the first equation in \for{critical} depends on the potential 
$\cV$. 

\Lem{eqrel-1} A nonzero vector $m=(\ml,\mr)\in(\rtre)^2$
satisfies $m\per\nabla_m\tilde H_\eta(b,m)=0$ iff either of the 
following alternative conditions is verified:
\bList
\item[(i)] $\ml$ and $\mr$ are both parallel to the same
principal axis; 
\item[(ii)] $\ml$ and $\mr$ both belong to the same principal
plane, say $e_i\oplus e_j$ with $b_i>b_j$, but are not colinear
with either principal axis, and, together with $b$, they satisfy
any of the following three conditions (where $(i,j,k)$ is a
permutation\footnote{By a permutation, we
mean either an even or an odd permutation.} of $(1,2,3)$):
\bListList
\item[(ii.1)] $b_j<|b_i-2b_k|$ and 
\begin{eqnarray}
   \Big({\ml^i-\mr^i\over \ml^i+\mr^i}\Big)^2 \ugarr 
   \Big({b_j+b_k\over b_j-b_k}\Big)^4 \, 
       {b_i^2-(b_j-2b_k)^2 \over b_i^2-(b_j+2b_k)^2 } 
\label{uv-1}  \\
   {\ml^j-\mr^j\over \ml^j+\mr^j} \ugarr
       -\, {\ml^i+\mr^i\over \ml^i-\mr^i}\,
       \Big({b_j+b_k\over b_j-b_k}\Big)^2 
       \Big({b_i+b_k\over b_i-b_k}\Big)^2 
       \,{b_i+b_j-2b_k\over b_i+b_j+2b_k}
        \;.
       \mathrm{\ \hbox to2truecm{}} 
\label{uv-2}
\end{eqnarray}
\item[(ii.2)] $b_i=2b_k+b_j$, $\ml^i=-\mr^i$, and $\ml^j=\mr^j$.
\item[(ii.3)] $b_i=2b_k-b_j$ and $\ml=\mr$.
\end{list}
\end{list}
\eLem

\acapo
Equations \for{uv-1} and \for{uv-2} 
imply that \for{uv-1} is also valid with the index $i$ and $j$ 
exchanged. If either $\ml$ or $\mr$ vanish (the case of an
irrotational ellipsoid), then equations \for{uv-1} and
\for{uv-2} are equivalent to $k=2$ and
$b_1^2b_2^2+b_2^2b_3^2+b_1^2b_3^2=3b_2^4$.

We now consider the equation $\nabla_bH_\eta=0$. Set
$$
   \tilde\cB_{j} \ug \big\{b\in\cB\,:\,\; 
           \MS_-(b_1,b_j,b_k)\ge 0 \,,\;
           \MS_+(b_1,b_j,b_k)\ge 0 \big\} \,,\qquad j=2,3
$$
where $(j,k)$ is any permutation of $(2,3)$ and 
$$
  \cB^\pm_{ij} \ug \big\{b\in\cB: 
            b_j\le \pm(b_i-2b_k) \,,\; 
            \D(b_i,b_j,b_k)\not=0 \,,\; 
            \MR_\mp(b_i,b_j,b_k)> 0 \,,\;
            \MR_\pm(b_j,b_i,b_k)>0 \big\}     
$$
where $(i,j,k)$ is any permutation of $(1,2,3)$ with $i<j$.
(We shall see below that some of these sets are empty.)

\Lem{eqrel-2} 
Assume $m=(\ml,\mr)\in(\rtre)^2$ is nonzero and satisfies
$m\per\nabla_m\tilde H_\eta=0$. Then $\nabla_b\tilde H_\eta=0$
iff one of the following two conditions is satisfied:
\bList
\item[(i)] $b\in\tilde\cB_{ij}$, $\ml=\ml^je_j$, $\mr=\mr^je_j$ 
and $(\ml^j\pm\mr^j)^2 = \MS_\pm(b_1,b_j,b_k)$
for either $j=2$, $k=3$ or $j=3$, $k=2$.
\item[(ii)] $b\in\cB^-_{ij}\cup\cB^+_{ij}$, 
$\ml=\ml^ie_i+\ml^je_j$, $\mr=\mr^ie_i+\mr^je_j$,
$(\ml^i\pm\mr^i)^2 = \MR_\pm(b_i,b_j,b_k)$, $(\ml^j\pm\mr^j)^2 =
\MR_\pm(b_j,b_i,b_k)$ and, if $b_i\not= 2b_k\pm b_j$, 
$$
    \ml^j-\mr^j = \sign\lp (2b_k-b_i-b_j) 
                {(\ml^j+\mr^j)(\ml^i+\mr^i)\over \ml^i-\mr^i} \rp \,
    \sqrt{\MR_-(b_j,b_i,b_k)} \,,
$$
for any permutation $(i,j,k)$ of $(1,2,3)$ with $i<j$.
\end{list}
\eLem

\acapo
The proof
of part (i) of Proposition~\ref{EqRel} is completed by showing
that the ``existence regions'' of the previous Lemma coincide
with those of the Proposition and that all the equilibria
described in that Lemma are \deq or adjoint to those indicated in
the Proposition:

\Lem{eqrel-3} The sets $\cB^+_{23}$, $\cB^-_{23}$ and
$\cB^-_{12}$ are empty, while $\tilde\cB_{2}= \cBSdue$,
$\tilde\cB_{3}= \cBStre$, $\cB^-_{13}= \cBI$,
$\cB^+_{13}=\cBII$, and
$\cB^+_{12}=\cBIII$.
\eLem

\Lem{eqrel-5} Under the assumptions of Lemma \ref{eqrel-2}, a
point $(b,0,m)$ is a relative equilibrium for some reduced
system iff it is \deq to one of the points 
$\big(b,0,(\mu_\a^\pm(b),\mu_\a^\mp(b))\big)$.
\eLem

\acapo
Part (ii) of the Proposition is a consequence of the following 

\Lem{eqrel-6}
The maps $\mu_\a^\pm$, $\a=\mathrm{S}_3,
\mathrm{II}, \mathrm{III}$, are never zero. For $\a=\mathrm{S}_2$
and $\a=\mathrm{I}$, they vanish on the sets indicated in
Proposition \ref{EqRel}.
\eLem

\acapo

\newsection{Appendix B: On Chandrasekhar's linearization}

To facilitate the comparison of our results with
Chandrasekhar's, in figure~\ref{whale} we display the regions of
elliptic ellipsoids of types II and III using the coordinates
employed by Chandrasekhar, namely $\big ({b_1 \over b_2}, {b_3
\over b_2} \big )$ for the ellipsoids of type II and $\big({b_1
\over b_3}, {b_2 \over b_3}\big)$ for the ellipsoids of type
III. (For the latter, Chandrasekhar also uses the coordinates
$\big({b_2 \over b_3}, {b_2 \over b_1}\big)$; the corresponding
figure should be compared to figure \ref{Ell-1}.b.) Among all
these ellipsoids, Chandrasekhar identified as elliptic
(``stable'', in his terminology) only the ellipsoids in the
lowest ellipticity region in figure~\ref{whale}.d.

\begin{figure}
\begin{picture}(6, 3.5)
\put(0.3,1.7){\scalebox{.7}{\includegraphics*{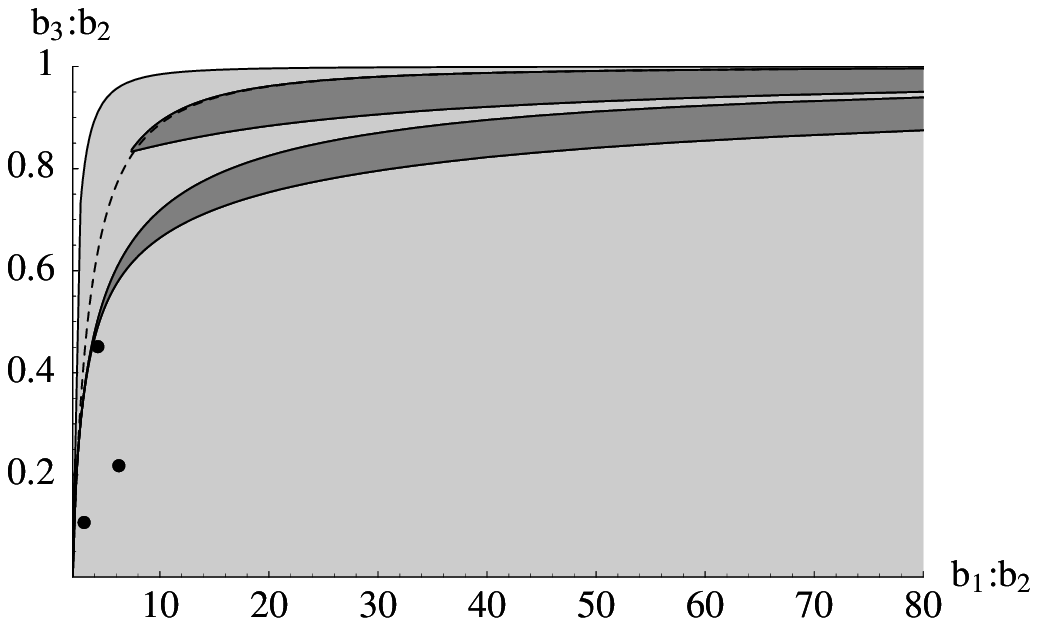}}}
\put(3.5,1.7){\scalebox{.7}{\includegraphics*{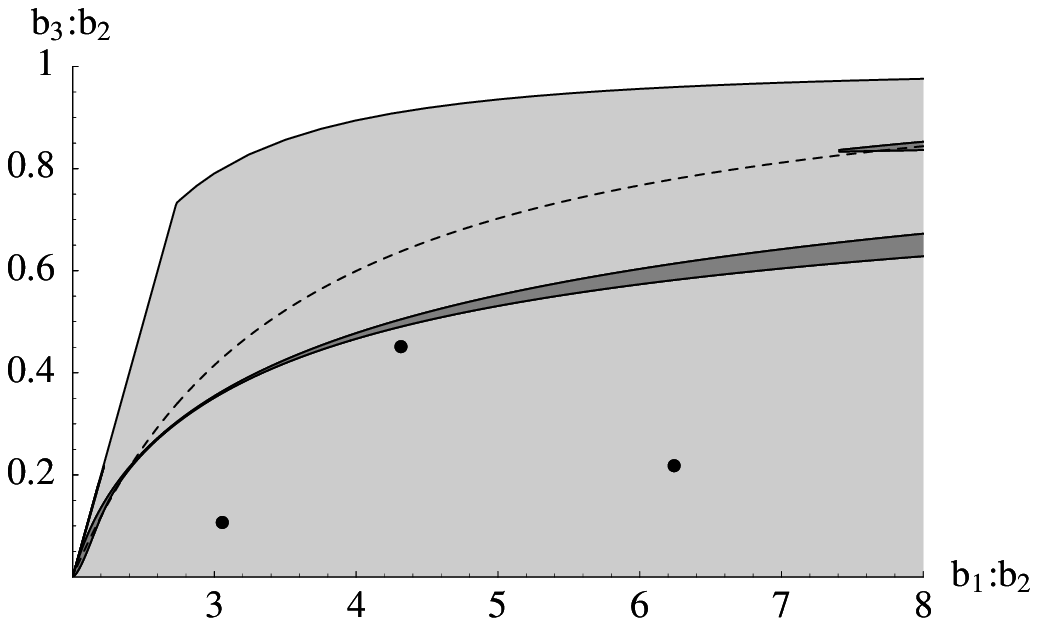}}}
\put(-0.25,0.2){\scalebox{.95}{\includegraphics*{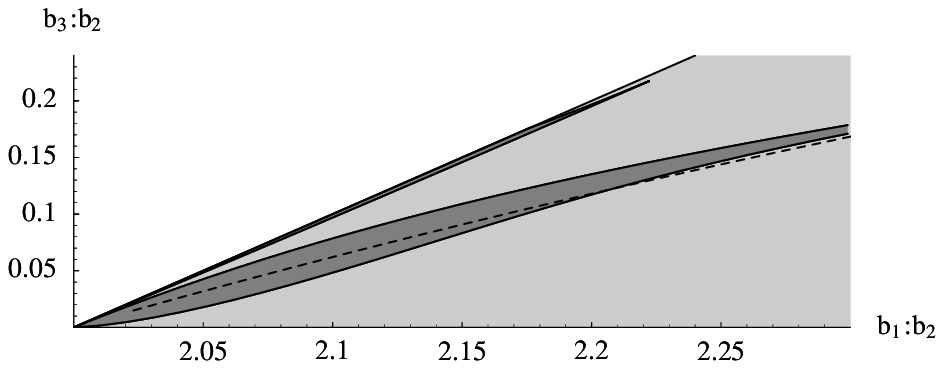}}}
\put(3.6,0){\scalebox{.75}{\includegraphics*{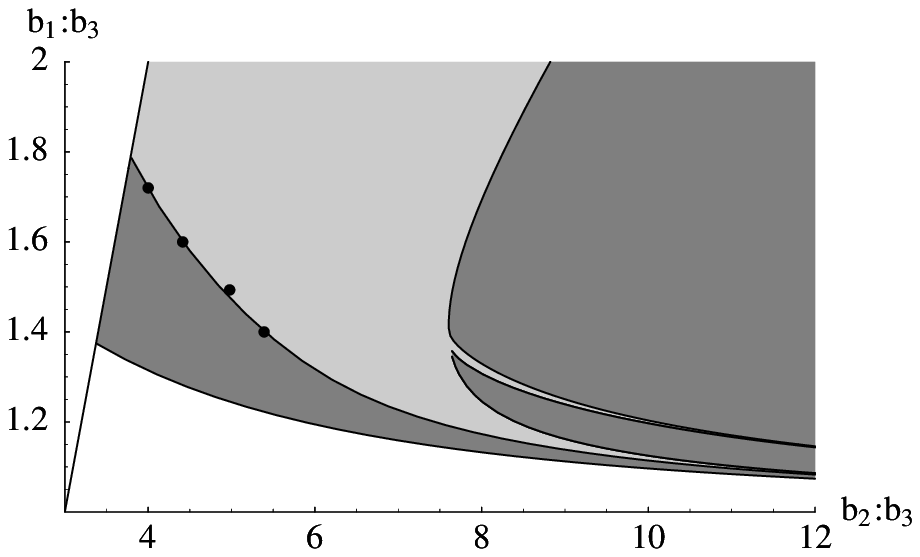}}}
\put(3.2,2.5){{\small (a)}}
\put(6.2,2.5){{\small (b)}}
\put(3.2,.8){{\small (c)}}
\put(6.2,.8){{\small (d)}}
\end{picture}
\caption{{\small Regions of elliptic (dark gray) and
non--elliptic (light gray) ellipsoids of types II (a--c), shown
at three different resolutions, and III (d). The indicated
individual points are the points tabulated in
\bib{chandrasekhar66}. The dashed curve in (a--c) is the zero
pressure curve.}}
\label{whale}
\end{figure}


As we have already remarked, Chandrasekhar's linearization of
Riemann's equations yields the same results as ours (after correcting
a few misprints that we report in the Remark below for the
convenience of the interested reader). Hence, it seems probable
that the meshes used by Chandrasekhar in his numerical study
were too coarse, or too small, to capture all of the ellipticity
regions shown in figure~\ref{whale}. Certainly, the very
few ellipsoids whose frequencies are reported in
\bib{chandrasekhar66} (three of type~II and four of type~III,
which are indicated by large points in figure~\ref{whale})
are consistent with our result.

Chandrasekhar did not exploit (nor, apparently, notice) the
Lagrangian and Hamiltonian structure of the Dirichlet problem
and worked with the full sixteen--dimensional Riemann's system,
without exploiting the symmetry of the problem. The
linearization of the reduced and unreduced systems are, of
course, equivalent: the four additional frequencies appearing in
the full system (two of which are identically zero) do not alter
the ellipticity.

\acapo
{\em Remark: } In \bib{chandrasekhar66}, the frequencies of the
linearized equations of motion are determined by the system of
nine equations (116). The nine-by-nine determinant of the
coefficients of this system is given in formula (133) of
\bib{chandrasekhar66}, where there are two misprints: the entry
(7,2) should be multiplied by $2+(\beta-2)\gamma {a_1^2/a_2^2}$,
and the entry (8,3) should be multiplied by
$2+(\gamma-2)\beta{a_1^2/a_3^2}$. This determinant is reduced to
the eight-by-eight determinant (134), in which the entry (1,8)
should be changed to $\O_2\O_3+2{a_1\over a_2}\O_2\O_3^\dagger$,
the entry (2,7) to $\O_2\O_3+2{a_1\over a_3}\O_2^\dagger\O_3$
and the entry (3,5) to $\O_2\O_3+{a_2\over
a_3}\O_2^\dagger\O_3^\dagger$; the entries (8,1), (7,2) and
(5,3) should be changed accordingly: the entry $(i,j)$ is
obtained from the entry $(j,i)$ by exchanging $\O_k$ and
$\O_k^\dagger$ ($k=2,3$). Also, the symbol $B_{123}$ appearing
in equation (135) of \bib{chandrasekhar66}, should be defined as
in equation (104), Chapter 3, of \bib{chandrasekhar69}. Let us
explicitly note that these misprints are not the sources of the
discrepancies between Chandrasekhar's and our conclusions: it is
only after they have been corrected that equation (134) of
\bib{chandrasekhar66} detects as elliptic the regions that
Chandrasekhar describes as ``stable''.

\vfill\eject


\end{document}